\documentclass[12pt,a4paper]{article}
\usepackage{graphicx,amsmath,amssymb} 
\usepackage{subfig,graphicx}
\usepackage{geometry}
\usepackage{tabularx}
\usepackage{rotating}
\usepackage{soul} % For \hl highlighting
\usepackage{xcolor}
\usepackage{lineno} % Add package in preamble
\usepackage[notcite,notref,color]{showkeys}

\usepackage{hyperref}

\usepackage{todonotes}

\newcommand{\R}{\mathbb{R}}
\newcommand{\id}{\mathbb{I}}

\newcommand{\VAR}{\mathrm{VAR}} % vectorautoregressive
\newcommand{\var}{\mathrm{Var}} % variance
\newcommand{\cov}{\mathrm{Cov}} %
\newcommand{\tr}{\mathrm{tr}}

\title{Early warnings of critical transitions through vector autoregression: lessons from multiscale systems}
\author{Bryony Hobden \and Paul Ritchie \and Peter Ashwin}
\date{March 2026}

\begin{document}

\maketitle

\begin{abstract}
In a nonautonomous nonlinear dynamical system, generic critical transitions (tipping points) are not limited to slow passage through fold bifurcations. They can also correspond to slow passage through other generic bifurcations, such as subcritical Hopf or even (in the multiscale setting) singular Hopf bifurcation. In systems such as these, generic early warning signals associated with critical slowing down, such as observing an increase in variance and lag-1 autocorrelation, can be obscured by oscillatory behaviour. Here, we present a method using vector autoregression (VAR) that aims to identify multiple eigenvalues of the system from multiple time series. With this, we are not only able to deduce the system's stability but also identify the type of bifurcation we are approaching and, therefore, the underlying dynamics, which gives us insight into the system's future state post-tipping. In practice, this can be beneficial for estimation and mitigation of risk. We look at three bifurcations in particular: the fold, the subcritical Hopf bifurcation, and the singular Hopf bifurcation in the presence of an extra time-scale separation. We hope to show that VAR is potentially a powerful tool in an EWS toolkit and can be more enlightening than simply looking at AR(1) of individual time series, especially if there is more than one timescale present in the system.
\end{abstract}

\newpage

\tableofcontents

\newpage

\section{Introduction}

Critical transitions or tipping points of nonlinear forced systems can be induced by slow passage through bifurcation points where an attracting equilibrium loses stability at some critical value of the forcing \cite{scheffer2001catastrophic}. Such bifurcation-induced transitions (B-tipping) can be distinguished from noise-induced (N-tipping) and rate-induced (R-tipping) \cite{ashwin2012tipping}. In this paper, we focus on generic mechanisms of bifurcation induced tipping for systems that may be multiscale.  

In the presence of noise, early warnings of loss of stability can be associated with the phenomenon of ``critical slowing down'' (CSD) (e.g. \cite{dakos2014critical,Dakos2008CSDClimate}), where an increase in variance and autocorrelation of a typical observable correlates to a slowly changing system approaching a bifurcation point. This is typically assumed to be a fold (saddle-node) bifurcation, although other bifurcations can lead to critical transitions \cite{Kuehn2013CSDVar,Williamson_2015,Bury_2021,huang2024deep}. 
One can think of CSD as a way to estimate the leading (i.e. least stable) eigenvalue of the Jacobian at a nearby equilibrium. This method has been widely applied to try and identify indicators of imminent B-tipping in a slowly changing system, to warn whether a system may or may not tip into a different state. It has been observed that, in some circumstances, the warning may be missing \cite{Hastings2010NoEWS,Morr2024,Ashwin2025}. Evaluating the skill of trends in CSD as early warnings of critical transitions can be challenging to extrapolate even in cases where CSD is applicable \cite{lehnertz2024time,Ashwin2025}, or may be misinterpreted in cases where there is N-tipping \cite{Ditlevsen2010}. Note that early warning signals can be useful, even if the result is negative: the important factor is that they are skilful so that false positives and false negatives can be avoided.  

There is no unique way to find an early warning signal from a time series, even when assumptions about timescale separation and noise are satisfied. One may need to detrend, choose a sampling rate and window size that has enough points to give a good estimation of stability while still remaining approximately stationarity, though note that fitting a nonlinear model can help avoid detrending \cite{Ditlevsen2023,kwasniok2024semiparametric}. For rapidly forced slow systems, CSD typically does not work, but other approaches may be possible \cite{huang2024deep,chapman2025quantifying,ritchie2026}. 

In most studies of bifurcation-induced tipping, only fold bifurcations have been considered. This occurs where a stable and unstable equilibrium collide and annihilate each other. This is often represented in diagrams showing 1D systems that can be written in terms of a gradient potential. However, at codimension 1 there are two generic bifurcations: fold (saddle-node) and Hopf bifurcations, and only the fold has such a gradient potential.  In Section \ref{sec:generic} we discuss which of these are considered ``safe'' and ``dangerous'' bifurcations in regards to tipping behaviour. It is very rare that Hopf bifurcations are considered for critical transitions, and very few studies consider early warnings for such bifurcations.

One study \cite{Bury_2021} looks at deep learning for early warning signals in a range of bifurcations including Hopf. The paper highlights that normal forms of systems can reveal information about how the system will behave past the bifurcation or tipping point, i.e. if it will shoot off and go to an alternative stable state or oscillate. However, generic early warning signals extracted from a time series will not reveal this type of information. Hence a deep learning algorithm is able to identify normal forms in a range of applications.

The paper \cite{Williamson_2015} considers using vector autoregression (VAR) to obtain indicators for the approach to various types of bifurcation, in particular periodic attractors approaching homoclinic and Hopf bifurcations. For the homoclinic there are no signs of traditional critical slowing down, but a decreasing frequency which can be used as an indicator of approaching bifurcation. For the Hopf however, there is no significant frequency change but significant early warning with respect to the increasing real part of the eigenvalue. We build on \cite{Williamson_2015} and we go one step further to consider Hopf bifurcation in a multiple time scale system, which can be a Singular Hopf bifurcation if the instability involves both fast and slow modes.

To our knowledge, no study has proposed an early warning for a system exhibiting a singular Hopf bifurcation. However, fairly simple models in climate applications such as climate can certainly demonstrate Hopf instabilities (e.g. \cite{Alkhayoun2019} in a 5-box model of Atlantic Meridional Overturning Circulation AMOC). Moreover, the singular Hopf bifurcation is already present in simple models (e.g. \cite{roberts2016mixedmode} for a model of ENSO). 

In previous work, \cite{hobden2025regularization} show that a regularization of the conceptual model \cite{boers2018ocean} for Dansgaard-Oeschger (DO) events occurs via passage through a supercritical singular Hopf bifurcation at the onset of a DO event. The study \cite{hobden2025regularization} highlights that because of the singular Hopf instability, the predictability of tipping in the model using generic EWS is limited from the leading eigenvalue.

In this paper, we focus on understanding and quantifying generic instabilities (using vector autoregressive processes) that may be associated with B-tipping from an attracting equilibrium in a multiscale system. Section~\ref{sec:generic} reviews some basic considerations from one-parameter bifurcation theory before highlighting that in multiple timescale systems, this needs interpretation due to the presence of eigenvalues with differing timescales. This means we can identify three distinct generic mechanisms where an equilibrium can lose stability. 

Section~\ref{sec:estimators} recalls the $\VAR(k,p)$ process and how $\VAR(k,1)$ can be used to estimate stability for an equilibrium in a $k$ dimensional system with noise. We consider not only the estimators but also the propagation of uncertainty through to the eigenvalues. Section~\ref{sec:examples} introduces a comparative example of each of the mechanisms discussed in Section~\ref{sec:generic} before applying the method of Section~\ref{sec:estimators} to stochastic versions of these examples in Section~\ref{sec:estimateexamples}. We finish with a brief discussion in Section~\ref{sec:discuss}.

\section{Generic instabilities in nonlinear systems}
\label{sec:generic}

We discuss the connection between codimension-one bifurcations \cite{Arnold1994,kuznetsov2023elements} of stable equilibria and tipping points in generic nonlinear systems in which a stable equilibrium loses stability. After highlighting that the only generic cases are fold and subcritical Hopf bifurcations, we recall that Hopf bifurcations occur in regular and singular forms in multiple-timescale systems.

\subsection{Instabilities of stable equilibria in nonlinear systems}

Consider a nonlinear system $x\in\R^n$
\begin{equation}
\dot{x}=f(x,\lambda)
\label{eq:ode}
\end{equation}
with a single parameter $\lambda\in\R$ but possibly many ($n\geq 1$) state variables, where $f$ is assumed smooth in both arguments. We say $X$ is an equilibrium for $\lambda$ of $f(X,\lambda)=0$ and $X$ is linearly stable if all eigenvalues of the Jacobian matrix $J=Df(X,\lambda)$ have negative real parts.

Suppose we have a family of linearly stable equilibria $X(\lambda)$ defined for $\lambda\in I$ and $I$ some open interval. The Implicit Function Theorem gives that this ``branch''  of equilibria continues smoothly on changing $\lambda$ \cite{kuznetsov2023elements}, while the Hartman-Grobman Theorem \cite{kuznetsov2023elements} gives that the equilibrium is locally stable (orbit equivalent to the linearized system) sufficiently close to the $X(\lambda)$. Moreover, the eigenvalues of the Jacobian matrix $Df(X(\lambda),\lambda)$ will change continuously with $\lambda$.

If $\lambda^*\in \partial I$ and
$$
x^*=\lim_{\lambda\rightarrow\lambda^*}X(\lambda)
$$
then $Df(x^*,\lambda^*)$ will necessarily have at least one ``critical'' eigenvalue with real part zero and we say the equilibrium $x^*$ is at bifurcation \cite{kuznetsov2023elements}. We say a bifurcation is {\em dangerous} \cite{Thompson2005} if there is a neighbourhood $N$ containing $x^*$ such that in any neighbourhood of $\lambda^*$ there are $\lambda$ with no attractors contained within $N$. Note that irreversibility is usually associated with the existence of a dangerous bifurcation, as it will be an exceptional case for any other attractor to become unstable exactly the same time as the branch ending at $x^*$ destabilizes.

For a system of the form (\ref{eq:ode}) but with a slowly drifting parameter, slow passage through a dangerous bifurcation gives rise to B-tipping as the system adjusts rapidly to a new family of attractors on loss of stability of the branch that was being tracked \cite{ashwin2012tipping}. Since it is a generic condition that a system with a critical eigenvalue has only one such eigenvalue, there are only two cases we need to consider when varying one parameter in a typical system of the form (\ref{eq:ode}). These are:
\begin{itemize}
\item {\bf Single real zero eigenvalue}\\ In this case, $x^*,\lambda^*$ locally undergoes a steady bifurcation. The centre manifold theorem can be applied to show that in some neighbourhood of the bifurcation, the dynamics can be described by a centre manifold that, in this case, is one-dimensional and locally attracting. Typically, the next highest order terms will be non-zero meaning the only generic case is a fold (saddle-node) bifurcation \cite{kuznetsov2023elements}, which is always dangerous and is the typical case considered for B-tipping.
\item {\bf Single imaginary pair of eigenvalues} In this case, $x^*,\lambda^*$ undergoes a Hopf bifurcation and the centre manifold theorem implies there is a locally attracting two dimensional invariant manifold on which (locally) either  (i) a stable limit cycle emerges after the bifurcation (supercritical Hopf), or (ii) an unstable limit cycle that forms the boundary of the basin of attraction of the equilibrium in the centre manifold shrinks to the equilibrium at the bifurcation leaving only an unstable equilibrium afterwards (subcritical bifurcation). Clearly, the Hopf bifurcation is dangerous only if it is subcritical \cite{KUEHN2011}. 
\end{itemize}
In summary, the dangerous bifurcations we need to consider as potentially causing B-tipping from a stable equilibrium, on varying one parameter (dangerous codimension one bifurcations), are typically only fold or subcritical Hopf bifurcations. 

Note that if two or more parameters are varied, then other cases can appear robustly - with $k$ parameters, one can find bifurcations of codimension up to $k$ at isolated points. Similarly, if there are extra constraints such as symmetries, other bifurcations are possible \cite{kuznetsov2023elements}. Nonetheless, if we choose a particular path through a high-dimensional parameter space, then we are still effectively at the one-parameter case.

\subsection{Instabilities of stable equilibria in multiple timescale systems}
\label{sec:multiscaleinstabilities}

Many models of physical systems, even in cases that can be described using low-dimensional dynamical models, involve interactions of several variables as in (\ref{eq:ode}), but moreover the processes, and hence the variables, may evolve on vastly differing timescales. 

We concentrate on an idealized setting where the dynamics are described by only two timescales, and the variables evolve on one of these two timescales. We say $x\in\R^{n_f}$ are {\em fast variables} and $y\in \R^{n_s}$ are {\em slow variables} so $(x,y)\in\R^{n}$ with $n=n_f+n_s$ and instead of (\ref{eq:ode}) we consider
\begin{equation}
\begin{aligned}
\epsilon \dot{x}&=f(x,y,\lambda)\\
\dot{y}&=g(x,y,\lambda)    .
\end{aligned}
\label{eq:odefastslow}
\end{equation}
The scalar quantity $\epsilon$ represents a ratio of timescales. We assume that $0<\epsilon\ll 1$ and as before a single parameter $\lambda\in\R$. There is an established theory of multiscale dynamics (see, e.g., \cite{kuehn2015multiple, KUEHN2011}), and these additional timescales can give rise to new phenomena and transitions by understanding the system as perturbed from the {\em singular} case with timescale ratio $\epsilon=0$. The critical manifold 
$$
\mathcal{C}=\{(x,y)~:~f(x,y)=0\}
$$
will be invariant in the limit $\epsilon=0$, and normally hyperbolic stable or unstable parts of this persist as invariant manifolds for sufficiently small $\epsilon$: we call the persisting stable part the ``slow manifold'', by Fenichel's theorem \cite{Fenichel1971, Fenichel1979}.

The Jacobian of (\ref{eq:odefastslow}) will have eigenvalues that are in two groups, a group of $n_f$ eigenvalues that are $O(\epsilon^{-1})$ associated with the fast variables, and one group of $n_s$ eigenvalues that are $O(1)$ associated with the slow variables. On the slow manifold, all these eigenvalues have negative real parts.

This means that for $0<\epsilon\ll 1$, generic codimension one instabilities of linearly stable equilibria of (\ref{eq:odefastslow}) can be of various types (and up to two dimensions, these are the only cases).
\begin{itemize}
\item {\bf Regular fold (on slow manifold)} where a single slow eigenvalue passes through zero on changing $\lambda$, any others are stable, and all fast eigenvalues remain strongly negative and of order $O(\epsilon^{-1})$ \cite{Arnold1994,kuznetsov2023elements}.

\item {\bf Regular subcritical Hopf (on slow manifold)} where a single imaginary pair of slow eigenvalues pass through zero on changing $\lambda$, any others are stable and all fast eigenvalues remain strongly negative and of order $O(\epsilon^{-1})$ \cite{Arnold1994} \cite{kuznetsov2023elements}.
\item {\bf Singular Hopf} this is the case where a single real fast eigenvalue passes through zero, combines with a single real slow eigenvalue to form a complex pair and passes through the imaginary axis. This is (confusingly) on the fold of the slow manifold, but is only associated with a fold of the fast subsystem, not the full system. This case is analysed in detail in \cite{Krupa2001}; see for example \cite{kuehn2015multiple,guckenheimer2012unfoldings}. A particular feature of this last bifurcation is that the complex pair appears for $|\lambda-\lambda^*|=O(\epsilon)$ and the associated oscillatory motion only appears in a neighbourhood of $x^*$ of size $O(\epsilon)$.
\end{itemize}

We illustrate examples of each of these instabilities in Section~\ref{sec:examples}.

\section{Estimation of eigenvalues using VAR}
\label{sec:estimators}

Assuming slow variation of a stable equilibrium driven by noise, one can estimate the stability of the system, and any trends, based on past and current information.
A $k$-dimensional vector autoregressive process with lag $p$, written as $\VAR(k,p)$, \cite{lutkepohl2013introduction} is a stochastic process $\{X_n\in\R^k\}$ satisfying
\begin{equation}
    X_{n+1} = c + A_1X_n + A_2 X_{n-1} +\dots +A_p X_{n-p+1} + \xi_n
    \label{eq:VARkp}
\end{equation}
for each $n$, where $c\in\R^k$ are constants, $A_{i}$ for $i=1,\ldots,p$ are $k\times k$ real matrices and $\xi_n\in\R^k$ are independent Gaussian random vectors with zero mean and variance $\Sigma$. In the examples that follow, we concentrate on the case $p=1$ so that (\ref{eq:VARkp}) reduces to the $\VAR(k,1)$ process
\begin{equation} 
\label{eq:VARk1}
    X_{n+1} = c + A X_n + \xi_n
\end{equation}
where $X_n\in\R^k$ and $A$ is a $k\times k$ matrix. This method includes the scalar AR(1) approach as the special case $\VAR(1,1)$ but as we see (and highlighted in \cite{Williamson_2015}) holds additional information about the instability that is missing from AR(1) estimations.

\subsection{Estimating equilibrium stability using VAR}

Suppose we have a nonautonomous SDE on $x\in\R^k$ with state- and time-dependent noise, of the form
\begin{equation}
\label{eq:SDEnonaut}
    dx= f(x,t)dt + g(x,t)dw(t)
\end{equation}
where $f:\R^k\times \R\rightarrow\R^k$ and $g:\R^k\times\R\rightarrow \R^{k\times k}$ is smooth and $w(t)$ is a vector of Gaussian white noise (Wiener) process with zero mean and covariance $\id_k$. 

If the system is autonomous, and the noise is additive and stationary, then we can write
\begin{equation}
\label{eq:SDE}
    dx= f(x)dt + dw(t)
\end{equation}
where $w(t)$ is a vector of Gaussian white noise (Wiener) process with zero mean and covariance $\Sigma$. 

Assume that $f(x,t)$ is slowly varying, so (\ref{eq:SDE}) is value and the state is close to an equilibrium $x^*\in\R^k$. Assume the noise is additive and small. Then as long as the sample path remains close to $x^*$ the solutions of (\ref{eq:SDE}) will be well-approximated by the linearization. Hence the fluctuations around the equilibrium are well-approximated by the Ornstein-Uhlenbeck process for $x\in\R^k$
\begin{equation}
\label{eq:SDElinear}
    dx = [Df(x^*) (x-x^*)]dt + dw
\end{equation}
where $Df(x^*)$ is the Jacobian of $f$ evaluated at $x^*$ and we assume $Df(x^*)$ is stable (all eigenvalues real part less than zero) so that there is a stationary measure for (\ref{eq:SDElinear}).

Note that the linear SDE (\ref{eq:SDElinear}) can be explicitly solved to give
\begin{equation}
X_{n+1}=\exp[Df(x^*)\Delta t ]X_n+ V_n
\end{equation}
where
$$
V_n=\int_{s=-\Delta t}^{0} \exp[Df(x^*)s]W((n+1)\Delta t+s)\,ds
$$
This means that
$A=\exp[Df(x^*)\Delta t]$ and so
\begin{equation}
    Df(x^*) = \frac{1}{\Delta t} \ln A
    \label{eq:DffromAln}
\end{equation}
where $\ln(A)$ is the matrix logarithm. For $\Delta t$ small enough, this will be the principal value and so well defined. Note in particular that if $B=A-\id_k$ is such that $\| B\|<1$ (which can be assured by taking $\Delta t$ small enough), where $\id_k$ is the $k\times k$ identity, then
\begin{equation}
Df(x^*)= \frac{1}{\Delta t} \left[ B - \frac{B^2}{2}+\frac{B^3}{3}-\cdots\right]
\label{eq:DffromB}
\end{equation}

Alternatively, if we time-discretise (\ref{eq:SDElinear}) using an Euler- Maruyama method \cite{gardiner2004handbook} with timestep $\Delta t$ then the fluctuations of the solution can be approximated (for small $\Delta t>0$) by $\{X_n\in\R^k\}$ where $X_n=x(n\Delta t)-x^*$ and
\begin{equation}
\label{eq:EMSDElinear}
    X_{n+1} = X_n+ Df(x^* )X_n \Delta t + \sqrt{\Delta t}W_n
\end{equation}
where $W_n$ are independent increments of $w(t+\Delta t)-w(t)$ and so are zero mean Gaussian with covariance $\Sigma$. 
Let us write
\begin{equation}
A=\id_k+ Df(x^*)\Delta t\approx \exp(Df(x^*)\Delta t)
\label{eq:AfromDf}
\end{equation}
and
%where $\id_k$ is the $k\times k$ identity, and
$$
\xi_n= \sqrt{\Delta t} W_n
$$
is a Gaussian process with zero mean and covariance $\Sigma \Delta t$.
Note that (\ref{eq:EMSDElinear}) corresponds to the $\VAR(k,1)$ process in (\ref{eq:VARk1}).

Note that we can invert (\ref{eq:AfromDf}) to give the approximation
\begin{equation}
Df(x^*)= \frac{A-\id_k}{\Delta t}
\label{eq:DffromA}
\end{equation}
which corresponds to the linear truncation of (\ref{eq:DffromB}).
Hence, if we have an estimator $\tilde{A}$ of $A$ from fitting the fluctuations to (\ref{eq:VARk1}) then we have an estimator 
$$
\tilde{J}=\frac{\tilde{A}-\id_k}{\Delta t}
$$
for $Df(x^*)$. 

Consequently, the eigenvalues $\{\tilde{\mu_i}\}$ of $\tilde{J}$ can be seen as estimators for the eigenvalues $\{\mu_i\}$ of $Df(x^*)$. However, these estimators are typically biased (because of the nonlinearity of the expressions for the eigenvalues) and there will be several sources of uncertainty in this estimation. 

Even in the stationary case (\ref{eq:SDE}) sources of uncertainty include the following: truncation error in (\ref{eq:SDElinear}), discretization error in (\ref{eq:EMSDElinear}), and estimation error from $\tilde{A}$. 
Truncation error can be minimized by ensuring that the sample path is close to the equilibrium such that a linearization is a good enough approximation and that higher order terms can be neglected. Using a small enough time step of the simulation $\Delta t $ will reduce error in the discretization step. Finally, estimation error of $\tilde{A}$ can be reduced by using a large number of time series data points in an observation window over which $A$ is estimated. Section \ref{sec:estimateexamples} briefly discusses sources of uncertainty in the non-autonomous case and ways in which these errors can be reduced.

For the remainder of this paper, we focus on the specific case of eigenvalue estimation and uncertainty quantification for multiscale systems with two coupled dynamical variables, i.e. $k=2$.

\subsection{Uncertainty estimation of eigenvalues}

We present two methods for calculating the standard error on the eigenvalues of the estimator $\tilde{A}$ for the case $k=2$. Since the relationship between the VAR parameters and eigenvalues is nonlinear, the standard errors of the parameter estimates cannot be directly translated into eigenvalue uncertainties.

Suppose we are given standard errors for an estimator $\tilde{A}$ of $A$ such that
\begin{equation}
\label{eqn:AtildeSEtilde}
\tilde{A}=
    \begin{pmatrix}
    \tilde{a} &\tilde{b} \\
    \tilde{c}& \tilde{d}
\end{pmatrix},~~
SE(\tilde{A})=
\begin{pmatrix}
    SE_{\tilde{a}} & SE_{\tilde{b}}\\
    SE_{\tilde{c}} & SE_{\tilde{d}}
\end{pmatrix}.
\end{equation}
We wish to understand how this uncertainty propagates through to the estimators for the Jacobian, $\tilde{J}$, and for the eigenvalues of the Jacobian, $\mu_i$. Note that (\ref{eq:DffromAln}) gives
\begin{equation}
\label{eqn:jac}
    \tilde{J} = \frac{1}{\Delta t}\ln(\tilde{A}) = :
    \begin{pmatrix}
        a&b\\
        c&d
    \end{pmatrix}.
\end{equation}

If the estimator $\tilde{A}$ is close to the identity, then to first order approximation the standard error on the entries of the Jacobian matrix is given by
\begin{equation}
    SE(\tilde{J}) = SE(\tilde{A)}/\Delta t = :
    \begin{pmatrix}
       SE_{a} & SE_{b} \\
    SE_{c} & SE_{d}
    \end{pmatrix}.
\end{equation}

To find the eigenvalues $\mu_i$ of the estimator $\tilde{J}$, we solve the characteristic equation given by
\begin{equation}
    \det(\tilde{J}-\mu\id)=0
\end{equation}
for eigenvalues $\mu_1, \mu_2$ and the identity matrix $\id$, where we assume $\Re(\mu_1)\geq \Re(\mu_2)$.
Hence
%The characteristic equation in the case of a $2\times2 $ Jacobian is
\begin{equation}
    \mu^2 - (a+d)\mu+ (ad-bc) = 0
\end{equation}
which is equivalent to
\begin{equation}
    \mu^2 - \tr(\tilde{J})\lambda + \det(\tilde{J}) = 0.
\end{equation}

The solution to the equation is given by the quadratic equation, such that the eigenvalues are a function of the variables $a,b,c$ and $d$.

\begin{equation}
\label{eqn:fabcd}
    \Re(\mu_{1,2}) = F_{1,2}(a,b,c,d) = \left\{
    \begin{array}{cl} 
    \frac{(a+d) \pm \sqrt{(a+d)^2-4(ad-bc)}}{2} & \mbox{Case }I\\
    \frac{(a+d)}{2} & \mbox{Case } II
    \end{array}\right.
\end{equation}
where Case $I$ corresponds to the real case with $(a+d)^2-4(ad-bc)>0$ (subscripts $1,2$ correspond to taking $+$ and $-$ respectively) and Case $II$ corresponds to the complex case with $(a+d)^2-4(ad-bc)\leq0$ (and subscripts $1,2$ are the same quantity),
as
\begin{equation}
\mu_{1,2} = \frac{\tr(\tilde{J}) \pm \sqrt{(\tr(\tilde{J}))^2-4\det(\tilde{J})}}{2}.
\end{equation}
We want to know
\begin{equation}
SE(\Re(\mu_{1,2}))  = 
\sqrt{\var(\Re(\mu_{1,2}))} = \sqrt{\var(F_{1,2}(a,b,c,d))}
\end{equation}

One way to approach this problem is the ``delta'' method based on a first-order Taylor series approximation \cite[Section 5.5.4]{Casella_Berger_2024}, \cite[Chapter~3]{TaylorJohnR2022}. The delta method approximates the standard errors of transformations of random variables using the first-order Taylor approximation: see Appendix~\ref{app:delta}. Since $\tilde{J}$ is an estimator, the entries of the matrix $\tilde{J}$ are random variables. 
In this case we can apply this to the function $F_{1,2}(a,b,c,d)$ as in equation (\ref{eqn:fabcd}) by computing the first-order partial derivatives (writing $\partial_a=\frac{\partial}{\partial a}$ etc) which in Case I gives
\begin{equation}
\left.
\begin{aligned}
    \partial_a F_{1,2} &= \frac{1}{2} \pm \frac{a-d}{2}(a^2+d^2-2ad+4bc)^{-1/2}\\
    \partial_b F_{1,2} &= \pm c(a^2+d^2-2ad+4bc)^{-1/2}\\
    \partial_c F_{1,2} &= \pm b(a^2+d^2-2ad+4bc)^{-1/2}\\
    \partial_d F_{1,2} &= \frac{1}{2} \mp \frac{a-d}{2}(a^2+d^2-2ad+4bc)^{-1/2}
\end{aligned}\right\}
\end{equation}
and in Case II gives
\begin{equation}
\left.
\begin{aligned}
    \partial_a F_{1,2} &= \frac{1}{2}\\
    \partial_b F_{1,2} &= 0\\
    \partial_c F_{1,2} &= 0\\
    \partial_d F_{1,2} &= \frac{1}{2}
\end{aligned}\right\}
\end{equation}
This means we can find the standard error on the eigenvalues $\mu_{1}$ and $\mu_{2}$, which can be written to first order as 
\begin{align}
    SE(\mu_{1,2}) &= SE(F_{1,2}(a,b,c,d))\nonumber\\
    & = \sqrt{ \left(\partial_a F_{1,2} SE_a\right)^2
    +\left(\partial_b F_{1,2} SE_b\right)^2
    +\left(\partial_c F_{1,2} SE_c\right)^2
    +\left(\partial_d F_{1,2} SE_d\right)^2},
\end{align}
If the $SE$ terms are not small then second- and higher- order terms can be included: a higher-order delta method and a multivariate delta method are presented in \cite{Casella_Berger_2024}.

The second method we present, and arguably a much simpler way to propagate the standard error on the entries of matrix $\tilde{A}$ through to eigenvalues of $\tilde{J}$, is via a Monte Carlo simulation. This is the method used to calculate the standard errors for eigenvalues shown in Section~\ref{sec:estimateexamples}. 

Assume the entries of matrix $\tilde{A}$ to be independent random variables and that each element can be distributed according to a normal distribution with mean equal to the corresponding element of the given estimator matrix, $\tilde{A}$, and a standard deviation equal to the corresponding element of the standard error matrix, $SE(\tilde{A})$ (for example, we choose the first component distributed as  $\mathcal{N}(\tilde{a}, SE_{\tilde{a}}^2)$
with $\tilde{a}$ and $SE_{\tilde{a}}$ as in equation (\ref{eqn:AtildeSEtilde})). We sample from these distributions to produce a new realisation of the estimator $\tilde{A}$. The Jacobian is calculated as before using equation \ref{eqn:jac}, and the eigenvalues of the Jacobian can then be obtained. If entries of the matrix $\tilde{A}$ are sampled a large number of times, then one can produce a large sample of eigenvalues from which a sample mean and sample standard deviation can be calculated. For a sufficiently large sample size, the standard error on the real parts of the eigenvalues calculated by Monte Carlo simulation will be approximately equal to that calculated using the first-order delta method. In Section \ref{sec:estimateexamples}, the standard errors shown in the figures are computed using Monte Carlo simulation with 1,000 samples. Given the large number of simulations, these estimates closely approximate those obtained from the first-order delta method. We therefore report the Monte Carlo standard error estimates, though this choice is largely arbitrary.

\section{Examples of multiscale instabilities}
\label{sec:examples}

\begin{figure}[h]
\subfloat[$\lambda=0.2$]{\includegraphics[width=0.30\textwidth]{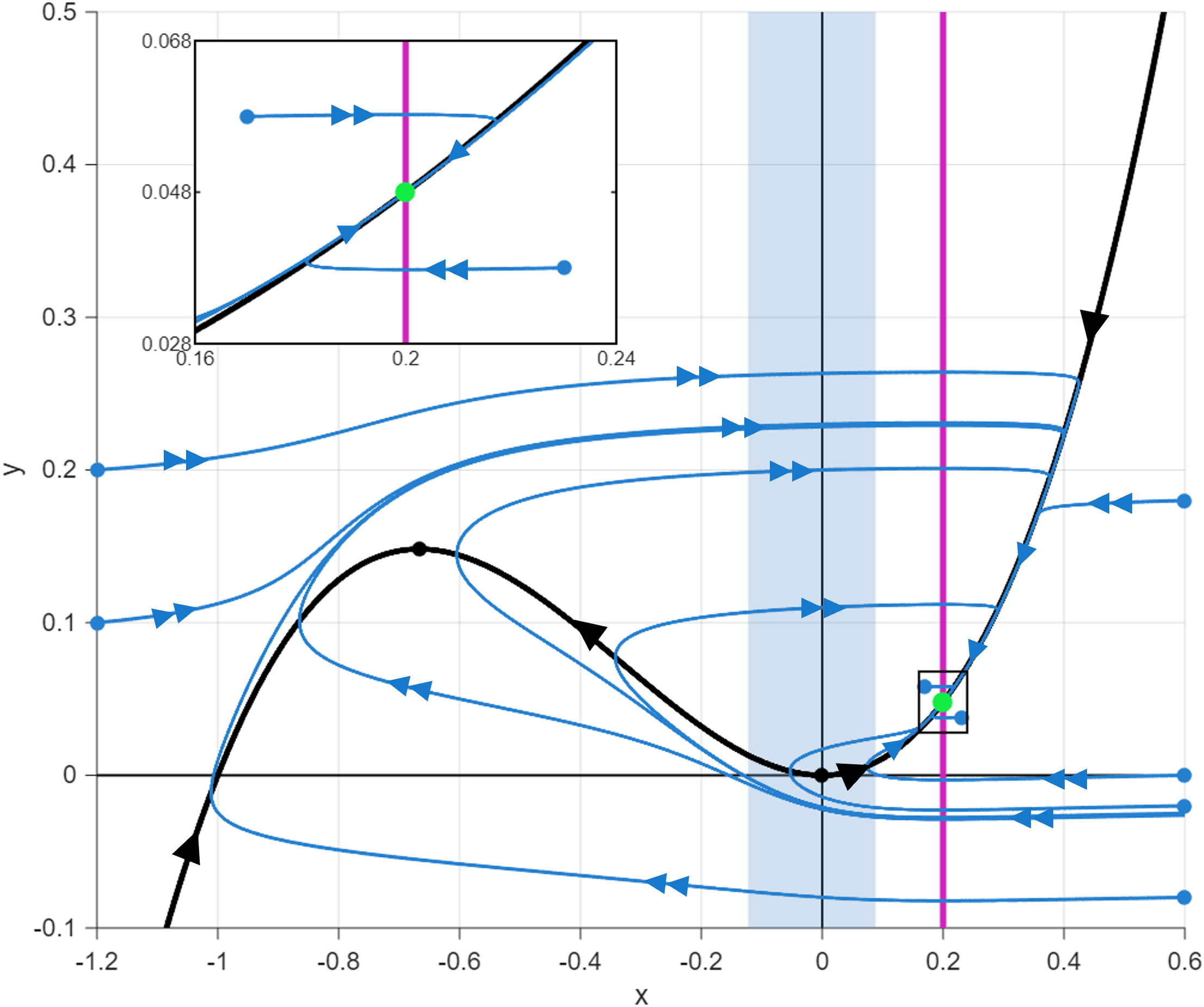}\label{fig:singhopfppa}}\hfill 
\subfloat[$\lambda=0.05$]{\includegraphics[width=0.30\textwidth]{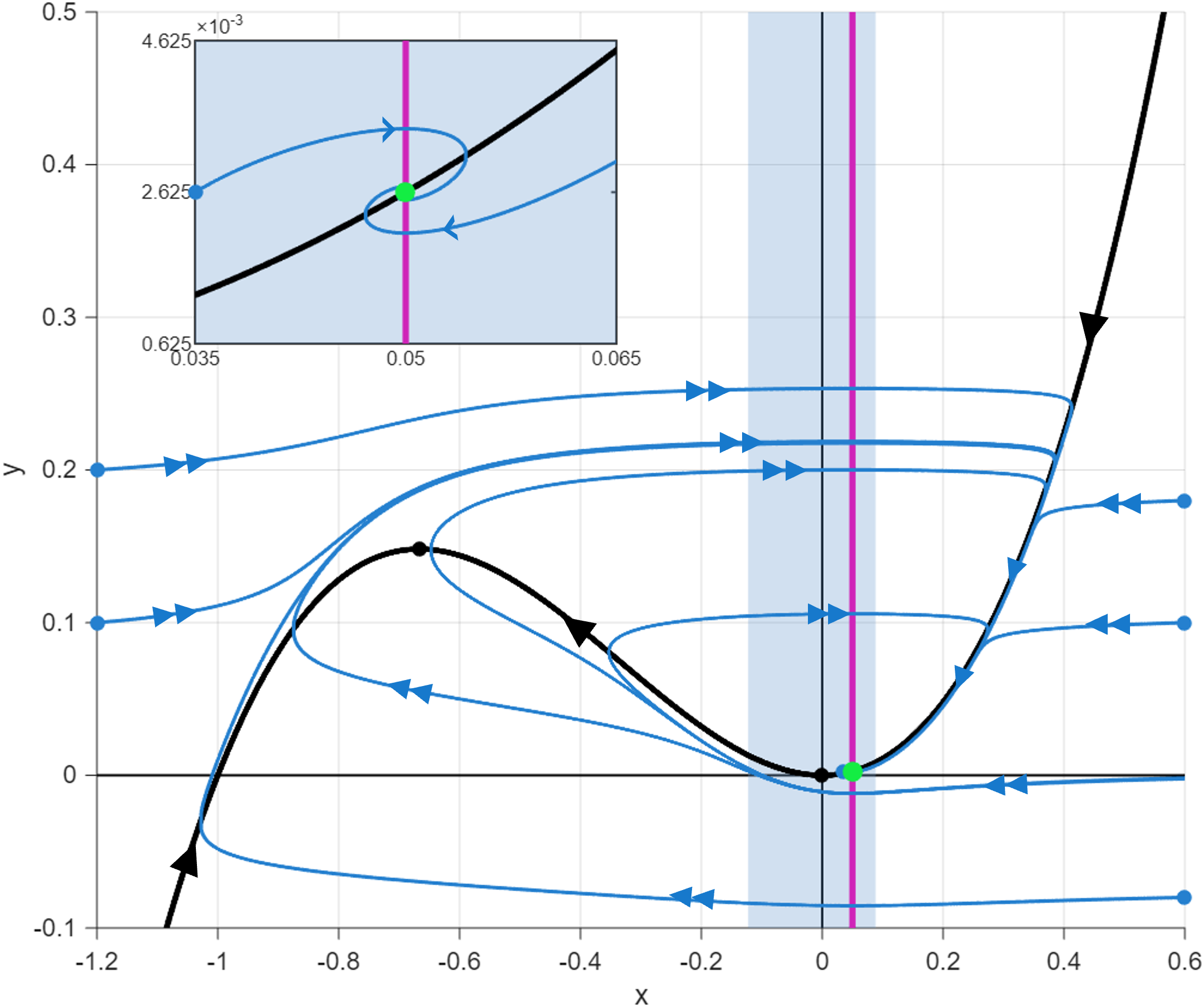}\label{fig:singhopfppb}}\hfill
\subfloat[$\lambda=0$]{\includegraphics[width=0.30\textwidth]{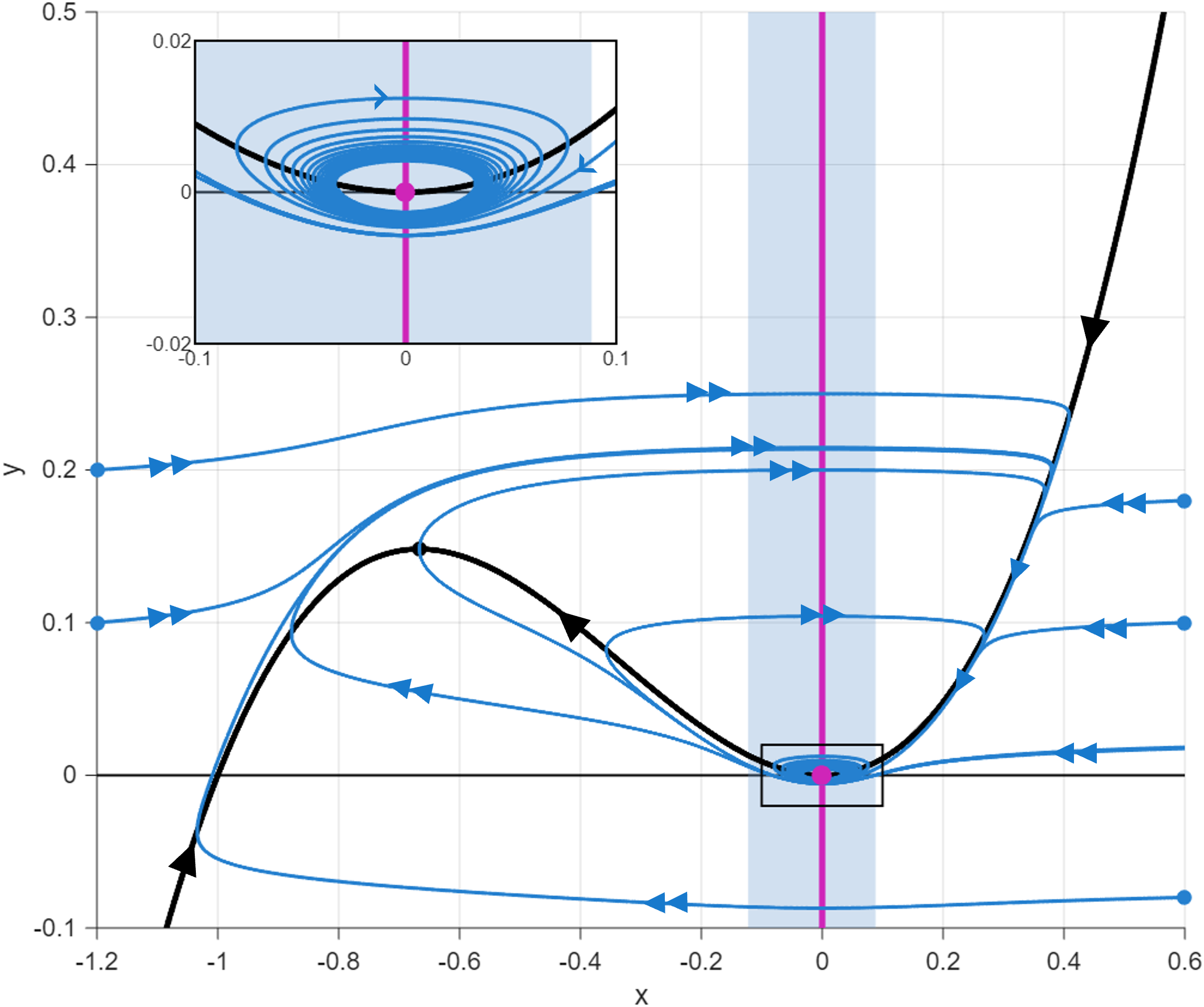}\label{fig:singhopfppc}}
    
\medskip
\hfill
\subfloat[$\lambda=-0.05$]{\includegraphics[width=0.30\textwidth]{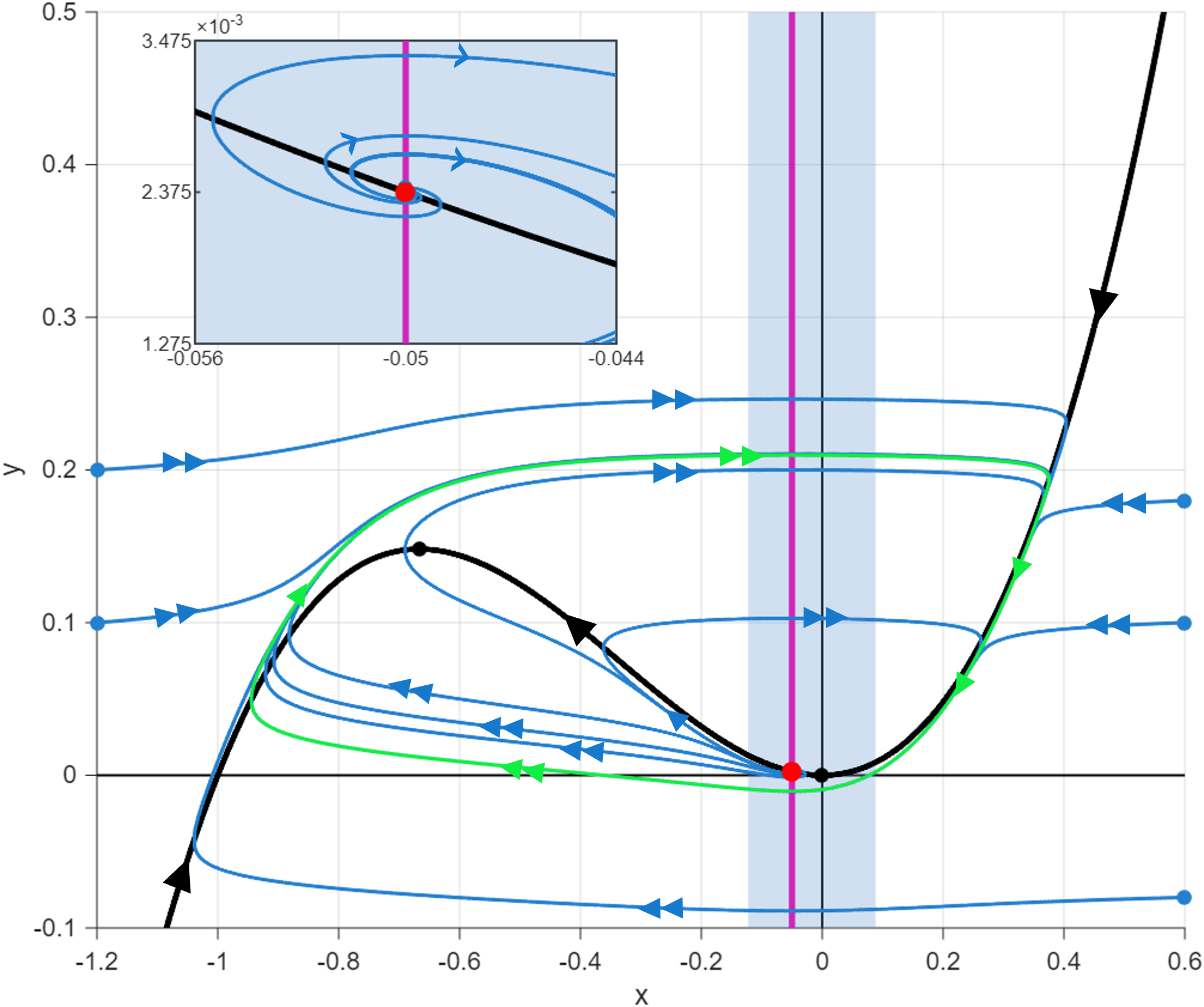}\label{fig:singhopfppd}}\hfill
\subfloat[$\lambda=-0.2$]{\includegraphics[width=0.30\textwidth]{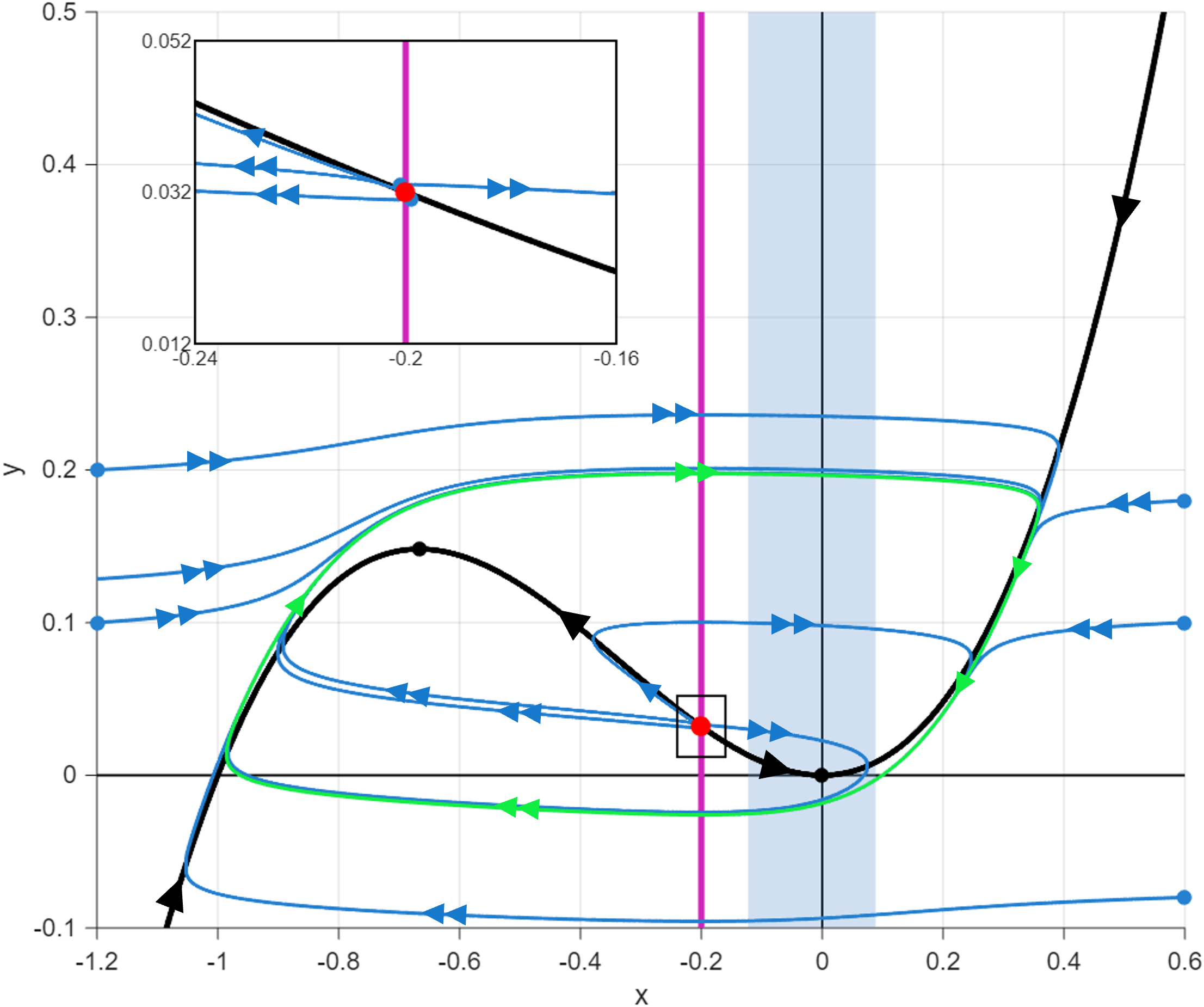}\label{fig:singhopfppe}}   
\hfill\mbox{}
\caption{Phase portrait diagrams for a singular Hopf bifurcation for the system (\ref{eq:singHopfxyslow}) for time scale separation $\epsilon=0.01$. The fast-slow system undergoes a singular Hopf bifurcation as the parameter $\lambda$ moves from right to left along the $ x$-axis. As the time-scale separation parameter $\epsilon$ decreases, the region where the system exhibits oscillatory behaviour shrinks (shaded in blue), scaling with $\sqrt{\epsilon}$. The thick black line is the critical manifold corresponding to $\dot{x}=0$, where single black arrows show the direction of travel along the critical manifold. Vertical magenta lines correspond to the nullcline for $\dot{y}=0$. Blue lines represent sample trajectories of the system, and blue points show the respective initial conditions. Arrows indicate the relative speeds of the trajectories. Insets show a zoomed-in panel around the equilibrium point for each $\lambda$. Subfigures \ref{fig:singhopfppb} and \ref{fig:singhopfppd}, the box around the region of interest is not shown due to scale. A stable limit cycle is shown in green in subfigures \ref{fig:singhopfppd} and \ref{fig:singhopfppe}.}
\label{fig:SingHopf}
\end{figure}

We now present examples of the three generic scenarios for instability outlined in Section~\ref{sec:multiscaleinstabilities}. In Section \ref{sec:estimateexamples} we assume that the forcing $\lambda(t)$ changes slowly with time. In Figures \ref{fig:Fold}, \ref{fig:Hopf} and \ref{fig:SingHopf}, we fix $\lambda$ and so the systems are autonomous. This allows us to analyse the equilibrium structure of each system for varying values of $\lambda$.

\subsection{Fold on slow manifold}

\begin{figure}[ht]
\subfloat[$\lambda=0.1$]{\includegraphics[width=0.31\textwidth]{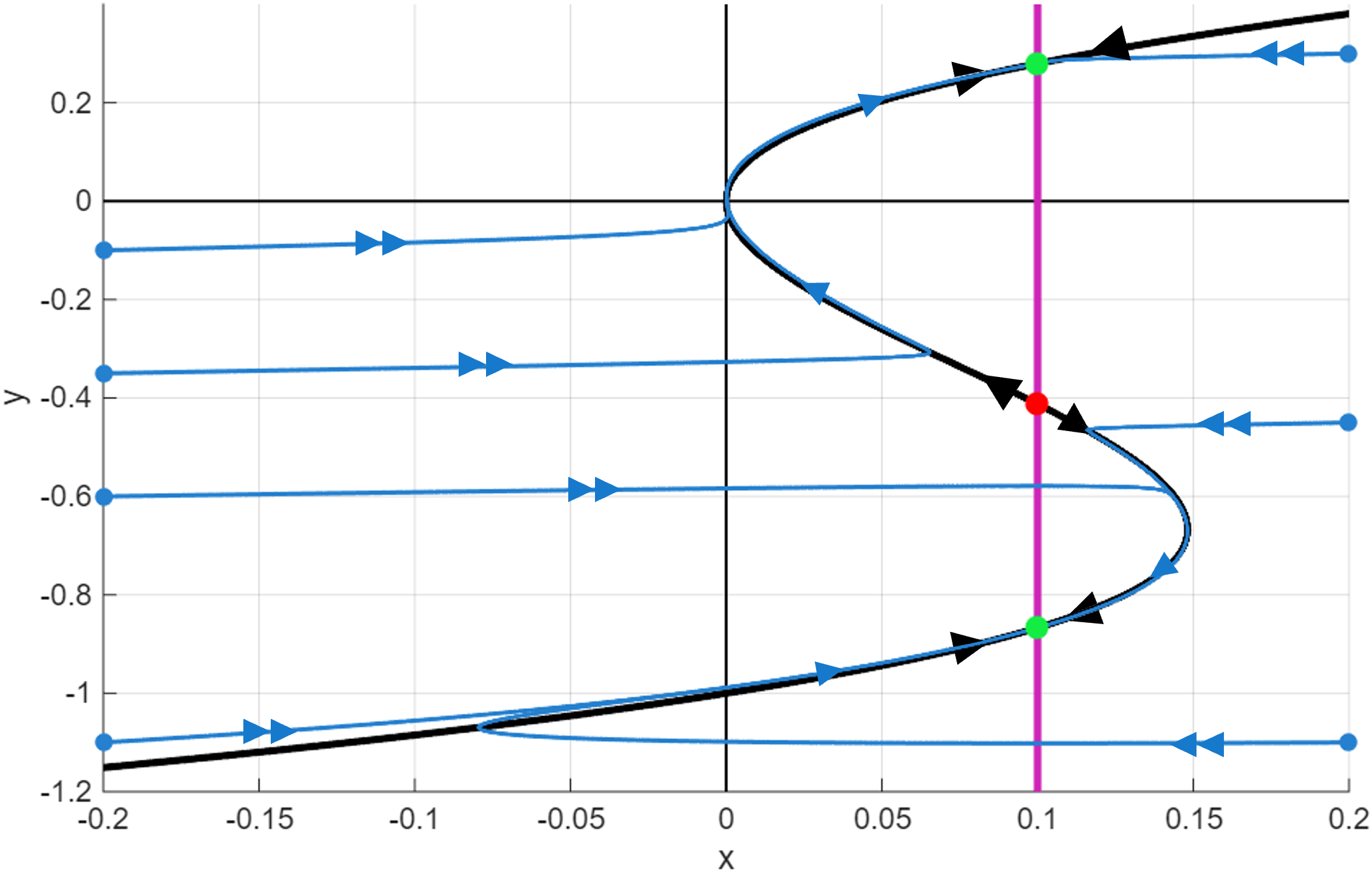}}\hfill 
\subfloat[$\lambda=0$]{\includegraphics[width=0.31\textwidth]{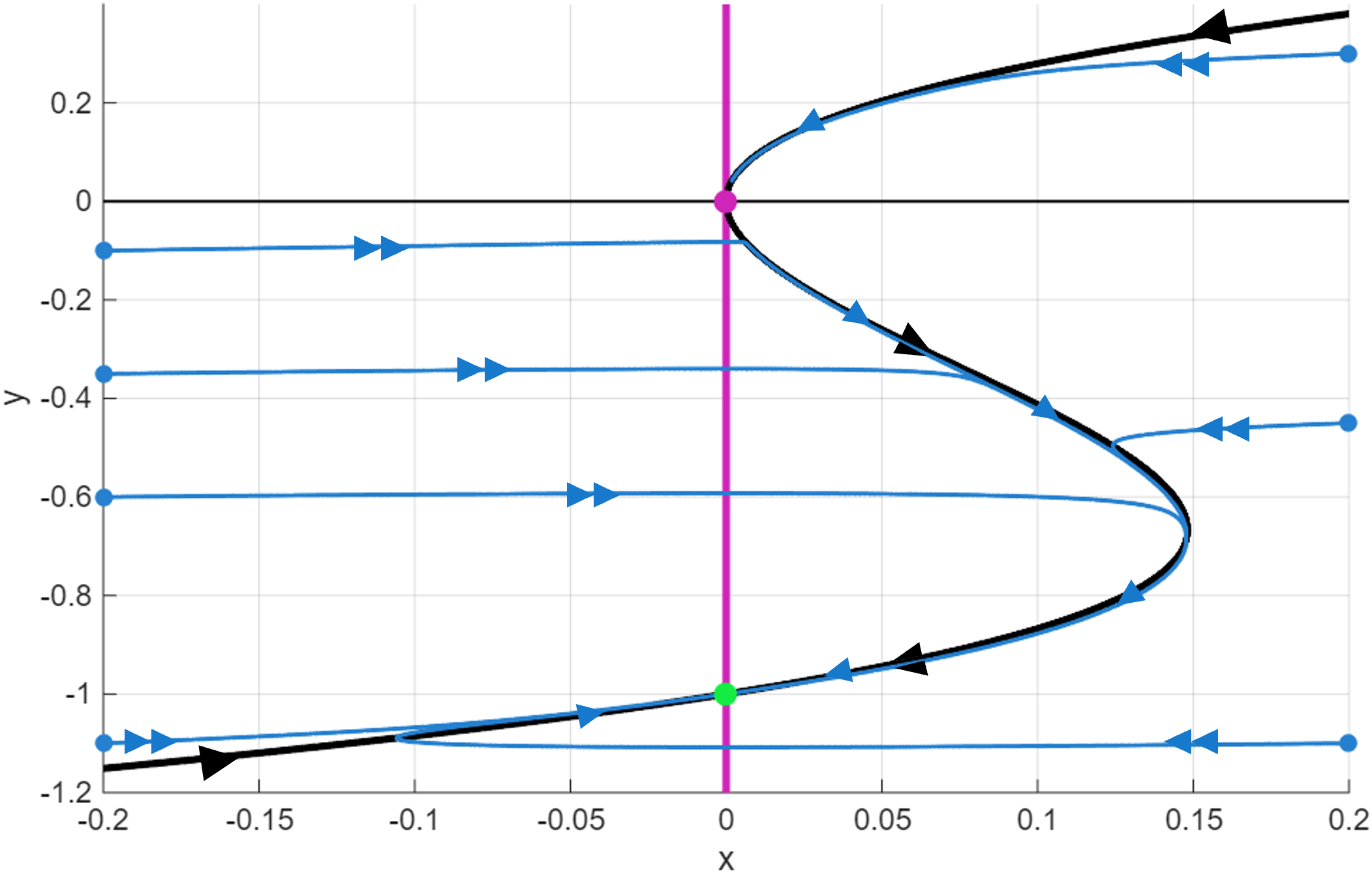}}\hfill
\subfloat[$\lambda=-0.1$]{\includegraphics[width=0.31\textwidth]{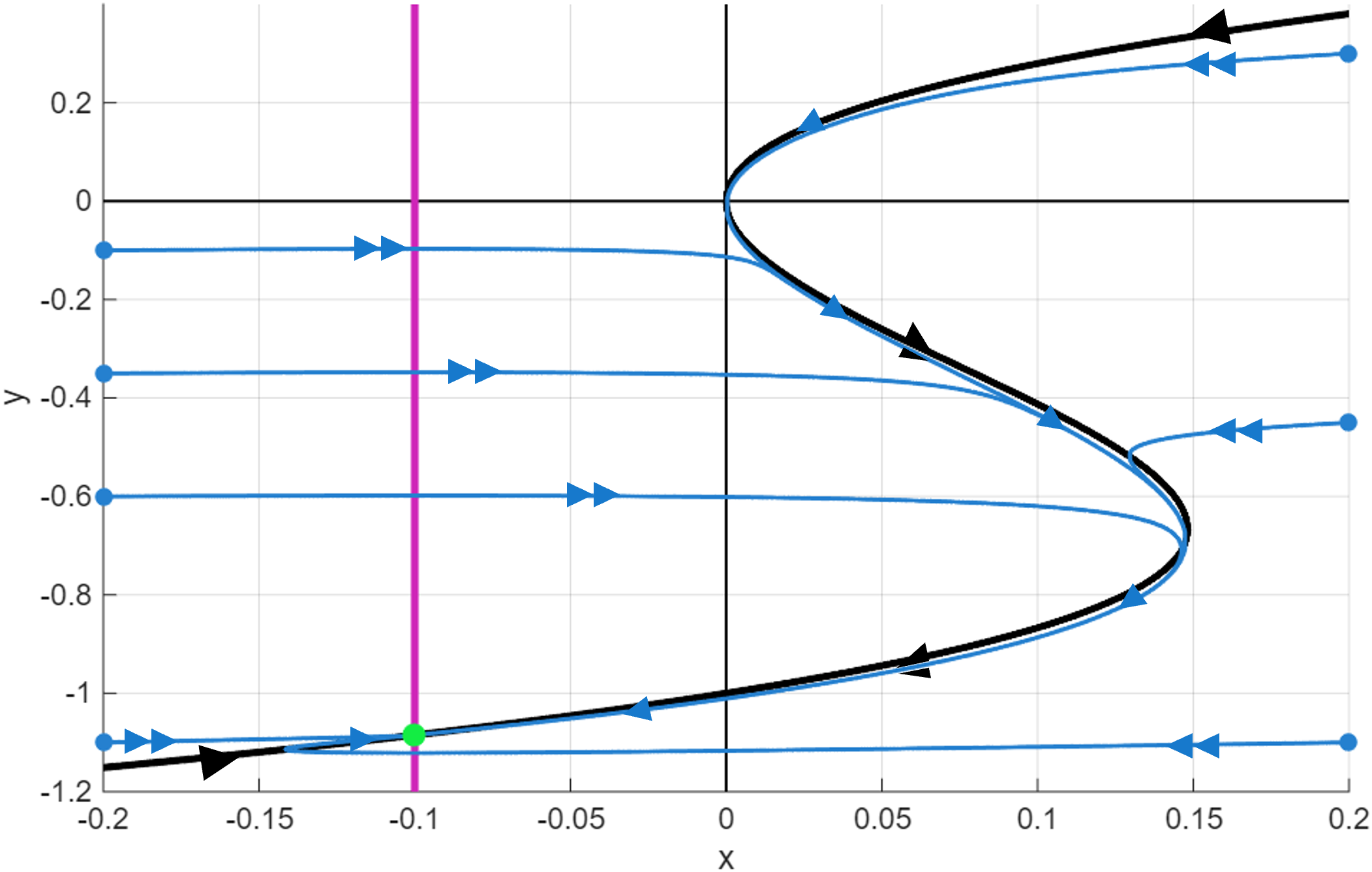}}
\caption{Phase portrait diagrams for a fold bifurcation for the system of equations (\ref{eq:fold_xy}) for time scale separation $\epsilon=0.1$. The fast-slow system undergoes a fold bifurcation at $(x,y)=(0,0)$ for $\lambda=0$. The thick black line is the critical manifold corresponding to $\dot{x}=0$, where single black arrows show the direction of travel along the critical manifold. Vertical magenta lines correspond to the nullcline for $\dot{y}=0$. Blue lines represent sample trajectories of the system, and blue points show the respective initial conditions. Arrows indicate the relative speeds of the trajectories.} \label{fig:Fold}
\end{figure}

An example of a multiscale system with a fold bifurcation on a (stable) slow manifold is given by the following set of equations:
\begin{equation}
\begin{aligned}
\label{eq:fold_xy}
    \epsilon \dot{x}&=y^2(1+y)-x\\
    \dot{y}&=\lambda-x,
\end{aligned}
\end{equation}
where $\epsilon\ll1$ is the timescale separation between the (fast) $x$ and (slow) $y$ variables. Note, we choose $r\ll 1 \ll \epsilon^{-1}$ so that the forcing timescale is much slower than both the $x$- and $y$-dynamics. Figure \ref{fig:Fold} shows phase portraits for various fixed values of $\lambda$ such that we see qualitatively different behaviour in the system. For small but positive values of $\lambda$, there are two stable equilibria separated by an unstable saddle point on the slow manifold. At the bifurcation point, $\lambda=0$, a stable equilibrium and the unstable saddle point collide with each other and annihilate to produce a non-hyperbolic point at $(x,y)=(0,0)$, so that only one stable equilibrium remains. For negative values of $\lambda$, the stable equilibria travel further in the negative $x$-direction along the slow manifold.

\subsection{Subcritical Hopf}

\begin{figure}[ht]
\subfloat[$\lambda=0.4$]{\includegraphics[width=0.33\textwidth]{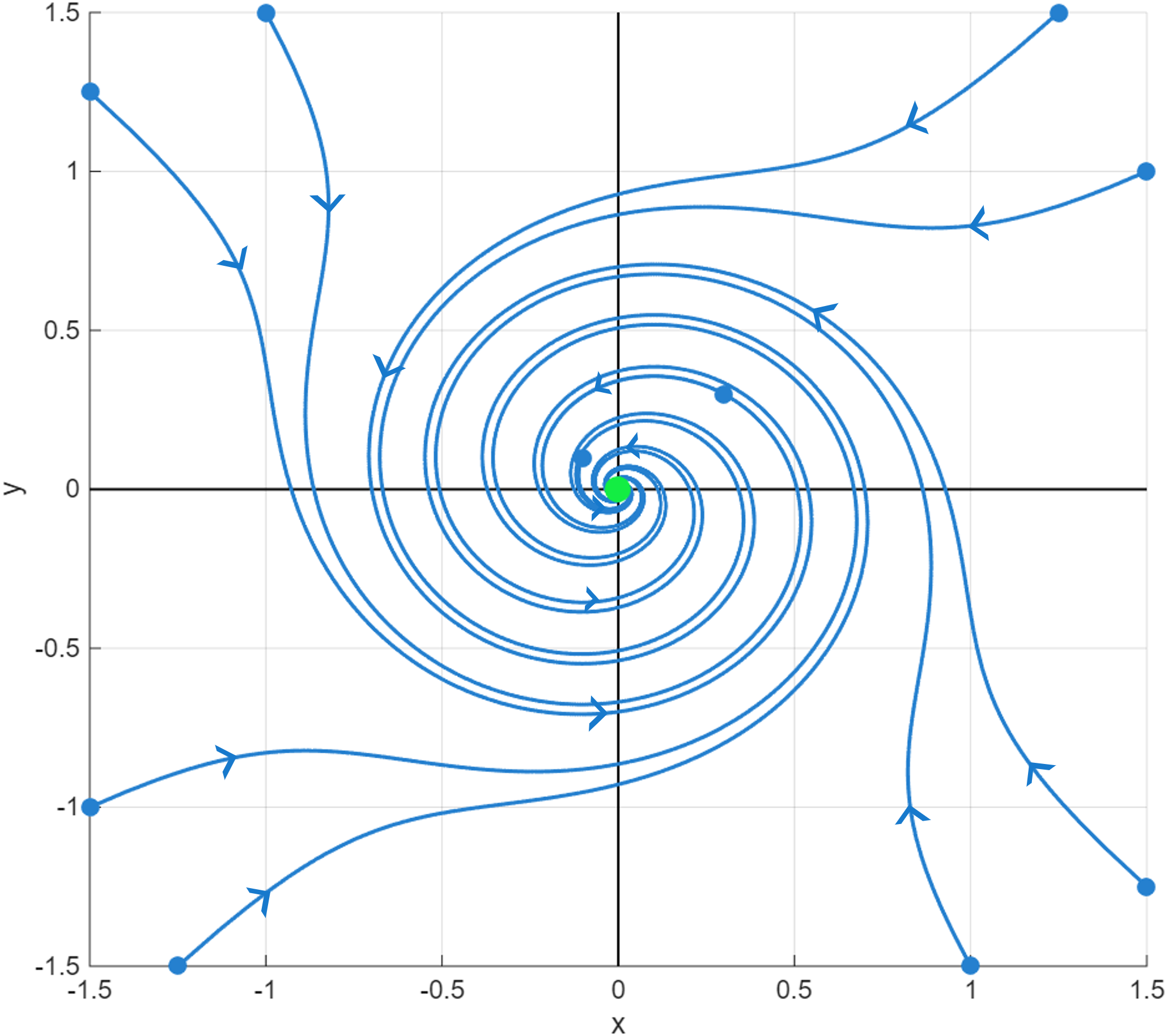}\label{fig:Hopfsubfiga}}
\subfloat[$\lambda=0.15$]{\includegraphics[width=0.33\textwidth]{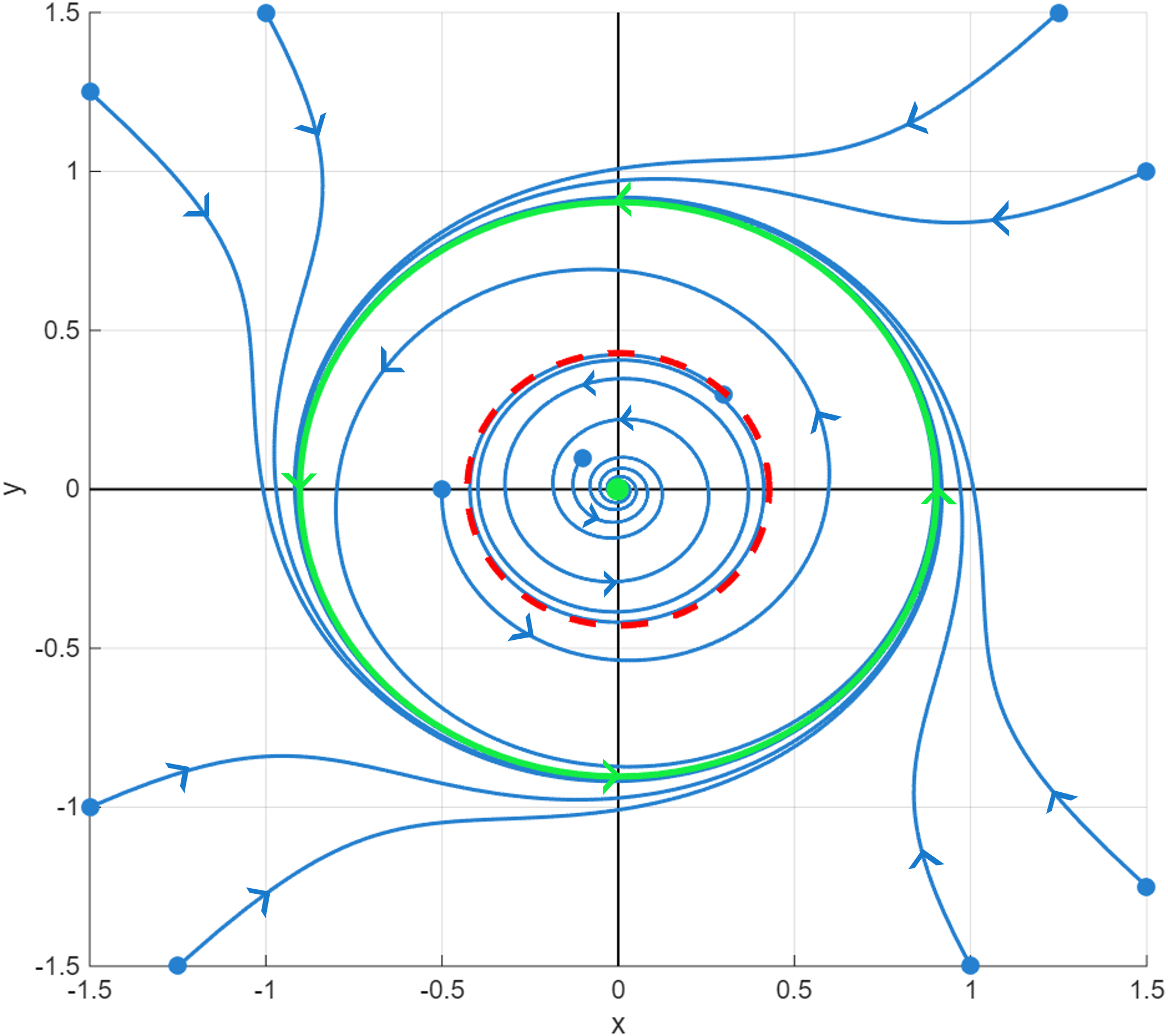}\label{fig:Hopfsubfigb}}\hfill
\subfloat[$\lambda=-0.25$]{\includegraphics[width=0.33\textwidth]{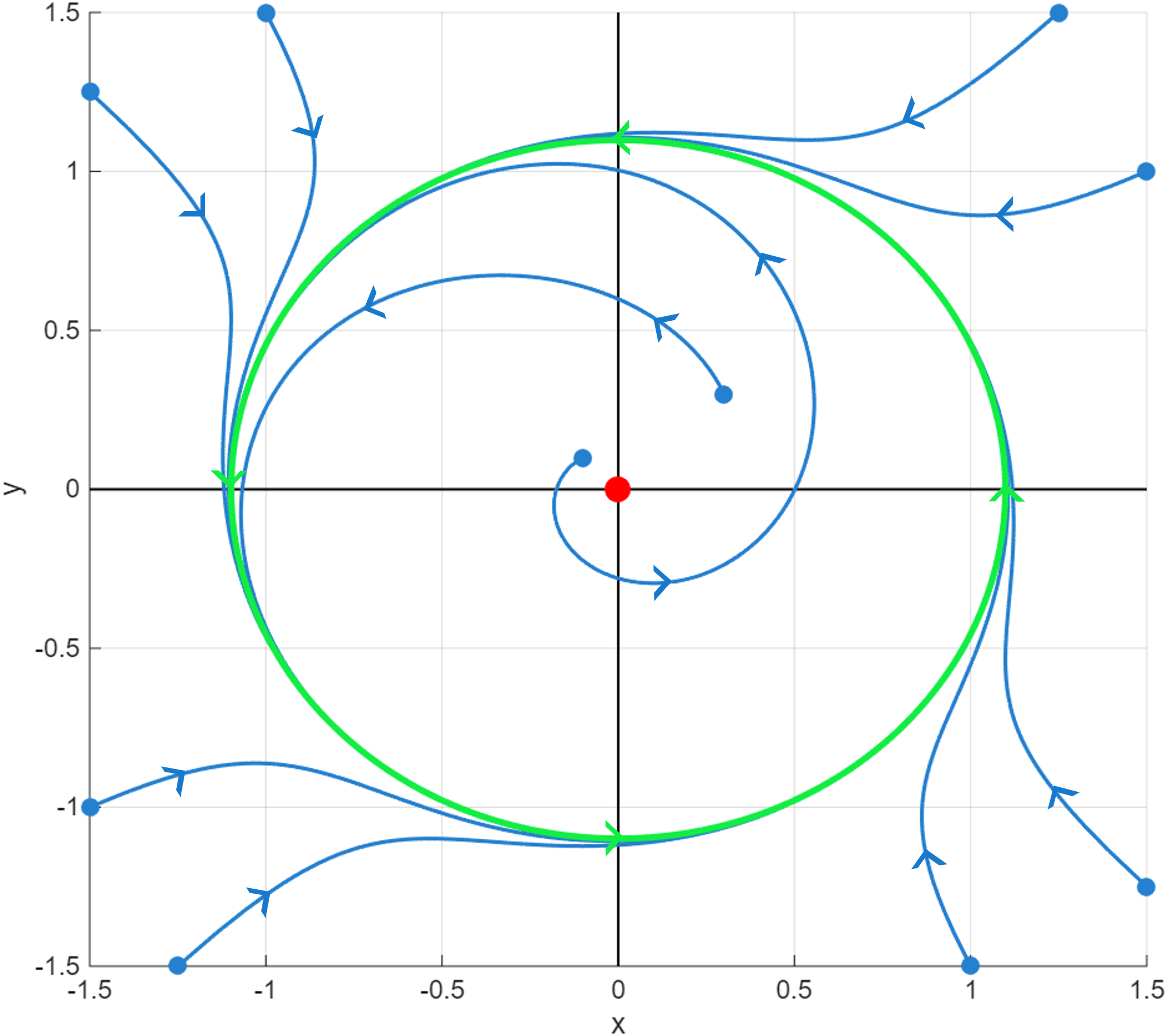}\label{fig:Hopfsubfigc}}\hfill 
\caption{Phase portrait diagrams for the subcritical Hopf bifurcation for the system in equation (\ref{eq:hopfz}). Blue lines represent sample trajectories of the system, and blue points show the respective initial conditions. In sub-figures (\ref{fig:Hopfsubfigb}) and (\ref{fig:Hopfsubfigc}), the green line indicates a stable limit cycle, and in sub-figure (\ref{fig:Hopfsubfigb}), the dashed red line indicates an unstable limit cycle. Green and red points represent stable and unstable nodes, respectively. For this system, the subcritical Hopf bifurcation occurs for $\lambda=0$. Here, the equilibrium at the origin becomes unstable, and the large amplitude limit cycle becomes the only local attractor. A saddle node bifurcation of periodic orbits occurs for $\lambda = 0.25$.}
\label{fig:Hopf}
\end{figure}

A system containing a subcritical Hopf bifurcation but nonetheless bounded dynamics is given by considering a path through a Bautin bifurcation normal form
\begin{equation}
\label{eq:hopfz}
    \dot{z}= (-\lambda+i\omega)z +|z|^2z - |z|^4z.
\end{equation}
In three or more dimensions, such a bifurcation may occur on a stable slow manifold where all transverse directions are fast and stable. We consider here only the simplest case where the slow manifold is the whole space. Writing equation (\ref{eq:hopfz}) in polar coordinates gives
\begin{align*}
\dot{\rho}&=-\lambda\rho +\rho^3 -\rho^5,\\
\dot{\theta}&=\omega
\end{align*}
which has explicit expressions for limit cycles at
\begin{equation}
    \rho_0=0, ~~ \rho_{\pm}=\sqrt{\frac{1\pm \sqrt{1-4\lambda}}{2}}.
\end{equation}
For $\lambda>1/4$, there are no real values of $\rho$ for the square root, and hence there is a single stable equilibrium at $\rho=\rho_0=0$. At $\lambda=1/4$, there is a single non-zero value for $\rho$, and here is where the saddle node bifurcation for periodic orbits occurs. The subcritical Hopf bifurcation (the switch of stability of the equilibrium point at the origin from stable to unstable) hence occurs at $\lambda=0$. For $0<\lambda<1/4$, the system is bistable, with a stable node at the origin and a large amplitude stable limit cycle $\rho=\rho_+$, separated by a small amplitude unstable limit cycle $\rho=\rho-$. In the ramping scenario, as the parameter $\lambda$ decreases, the system undergoes a dangerous bifurcation: the stability of the equilibrium at the origin is lost, and the system jumps to large-amplitude oscillations. This scenario of changing $\lambda$ is illustrated in phase portraits in Figure~\ref{fig:Hopf}.

\subsection{Singular Hopf}

Consider the following set of ordinary differential equations that, for a certain value of $\lambda$, will undergo a singular Hopf bifurcation given by
\begin{equation}
\begin{aligned}
\label{eq:singHopfxyfast}
     \epsilon\dot{x}&=y-x^2(1+x)\\ 
    \dot{y}&=\lambda-x. 
\end{aligned}
\end{equation}

It is clear that \ref{eq:singHopfxyfast} has an equilibrium point at $(x^*, y^*) = (\lambda, \lambda^2(1+\lambda))$ with a Hopf bifurcation occurring at $\lambda=0$. The Jacobian of the system is given by

\begin{equation}
J = 
    \begin{pmatrix}
        \frac{-2x-3x^2}{\epsilon} &\frac{1}{\epsilon}\\
        -1 &0
    \end{pmatrix}.
\end{equation}
At the Hopf bifurcation ($\lambda=0$) the eigenvalues of the Jacobian $J$ are $\mu_{1,2} = \pm i\sqrt{\frac{1}{\epsilon}}$ and $\mu_{1,2}$ will tend to infinity as the timescale separation parameter $\epsilon\rightarrow 0$. 
If we look at the system in terms of the fast time scale, i.e. $\tau = t/\epsilon $ such that $x' = \frac{dx}{d\tau}$, we have the following system of equations 
\begin{equation}
\begin{aligned}
\label{eq:singHopfxyslow}
    x'&=y-x^2(1+x)\\ 
    y'&=\epsilon(\lambda-x). 
\end{aligned}
\end{equation}

Looking at the system from the perspective of the fast timescale, the linearisation now produces eigenvalues $\mu_{1,2}=\pm i\sqrt{\epsilon}$, for which $\mu_{1,2}\rightarrow$ as $\epsilon \rightarrow0$. In each case, for both timescales, the eigenvalues at the Hopf bifurcation are singular, and this is how we can define a singular Hopf bifurcation. A more generic normal form of the Singular Hopf bifurcation in $\mathbb{R}^2$ is given in \cite[p.~205]{kuehn2015multiple}, but here we focus on a reduced form given in \ref{eq:singHopfxyslow}. 

For $\epsilon=1$, the system undergoes a supercritical Hopf bifurcation. In this section, we specifically look at the case $0<\epsilon\ll1$, which produces a singular Hopf bifurcation. Note that, unlike a Hopf bifurcation without timescale separation, it does not matter whether the singular Hopf bifurcation is subcritical or supercritical. Due to the large timescale separation in the singular Hopf case, and therefore to rapid movement away from the original stable state toward a large-amplitude periodic cycle, this can still be considered tipping behaviour. 

\section{A VAR approach to early warnings of tipping}
\label{sec:estimateexamples}
For each of our examples, we now assume that the forcing $\lambda(t)$ changes slowly with time, and is given by 
\begin{equation}
	\label{eq:ramp}
    \lambda(t) = \lambda_0 -rt,
\end{equation}
where $\lambda_0$ is some initial value and $r$ is a constant ramping rate with $0<r\ll 1$.

We use the method outlined in Section~\ref{sec:estimators} to give estimates of the eigenvalues for stochastic versions of the examples considered in Section~\ref{sec:examples} for the three bifurcation types. The parameters are chosen as in Table~\ref{tab:params} using the following considerations for each system:
\begin{itemize}
\item The rate parameter $r$ is chosen small enough to ensure slow ramping of the parameter, namely so that the system tracks the stable equilibrium up to the bifurcation point.
\item The noise on $x$ and $y$ variables is set to be independent, identical and sufficiently small, $\alpha_x = \alpha_y \ll 1$ that a noise-induced transition is not seen before passing the bifurcation point. 
\item The time series is subsampled to ensure stability of the method.
\item The block (window) size is chosen as a trade-off between (a) having enough data points to give an estimator with reasonably low standard error and (b) not having so many points that the nonstationarity within the window is minimized.
\end{itemize}
Table \ref{tab:params} summarises the parameter values used in the numerical simulations and for the VAR analysis, including subsampling and block size. We solve the system of stochastic differential equations using the Improved Euler method (Heun) \cite[p.~328]{suli2003introduction}. Alternative approaches, such as Strang splitting for higher-order stochastic systems, have been proposed in the literature \cite{pilipovic2025strang}.

\begin{table}
\centering
%\textcolor{red}{
\begin{tabularx}{\textwidth}{|c|X|c|c|c|}
\hline
\textbf{Parameter} & \textbf{Meaning} & \textbf{Fold} & \textbf{Subcritical Hopf} & \textbf{Singular Hopf} \\
\hline\hline
$\epsilon$ & Time scale separation & 0.1 & (1) & 0.01 \\
\hline
$r$ & Ramping rate & 0.001 & 0.01 & 0.005 \\
\hline
$\Delta t$ & Simulation time step & 0.002 & 0.01 & 0.001 \\
\hline
$\lambda_0$ & Initial value of $\lambda$ & 0.3 & 3.0 & 0.4 \\
\hline
$\alpha_x$ & Noise amplitude on the $x$ variable & 0.01 & 0.05 & 0.005 \\
\hline
$\alpha_y$ & Noise amplitude on the $y$ variable & 0.01 & 0.05 & 0.005 \\
\hline
$\omega$ & Rotational frequency & - & 0.3 & - \\
\hline
Tspan & Length of simulation & 350 & 350 & 90 \\
\hline
Sub & Sampling rate & 10 & 10 & 5 \\
\hline
Block Size & Number of points in observation window & 1250 & 250 & 1000 \\
\hline
\end{tabularx}
%}
\caption{Table of parameters used for the stochastic simulation of the fold, subcritical Hopf and singular Hopf systems as illustrated in Figures~\ref{fig:foldevs}, \ref{fig:hopfevs} and \ref{fig:singhopfevs}.}
\label{tab:params}
\end{table}

\subsection{VAR for fold on slow manifold}

To evaluate the stability of the nonautonomous fold system (\ref{eq:fold_xy}) using VAR techniques, we consider the following It\^{o} SDE with additive noise:
\begin{equation}
	\begin{aligned}
		\label{eq:fold_dxdy}
		\epsilon dx &=\big(y^2(1+y)-x\big)dt +\alpha_xdW_x\\
		dy&= (\lambda - x)dt + \alpha_ydW_y\\
		d\lambda &= -rdt
	\end{aligned}
\end{equation}
where $dW_x$ and $dW_y$ represent independent Wiener (white noise) processes with amplitudes $\alpha_x$ and $\alpha_y$ chosen sufficiently small such that transitions in the system are primarily due to the bifurcation and not a noise-induced transition. 
\begin{figure}[ht]
\centering
\includegraphics[width=\textwidth]{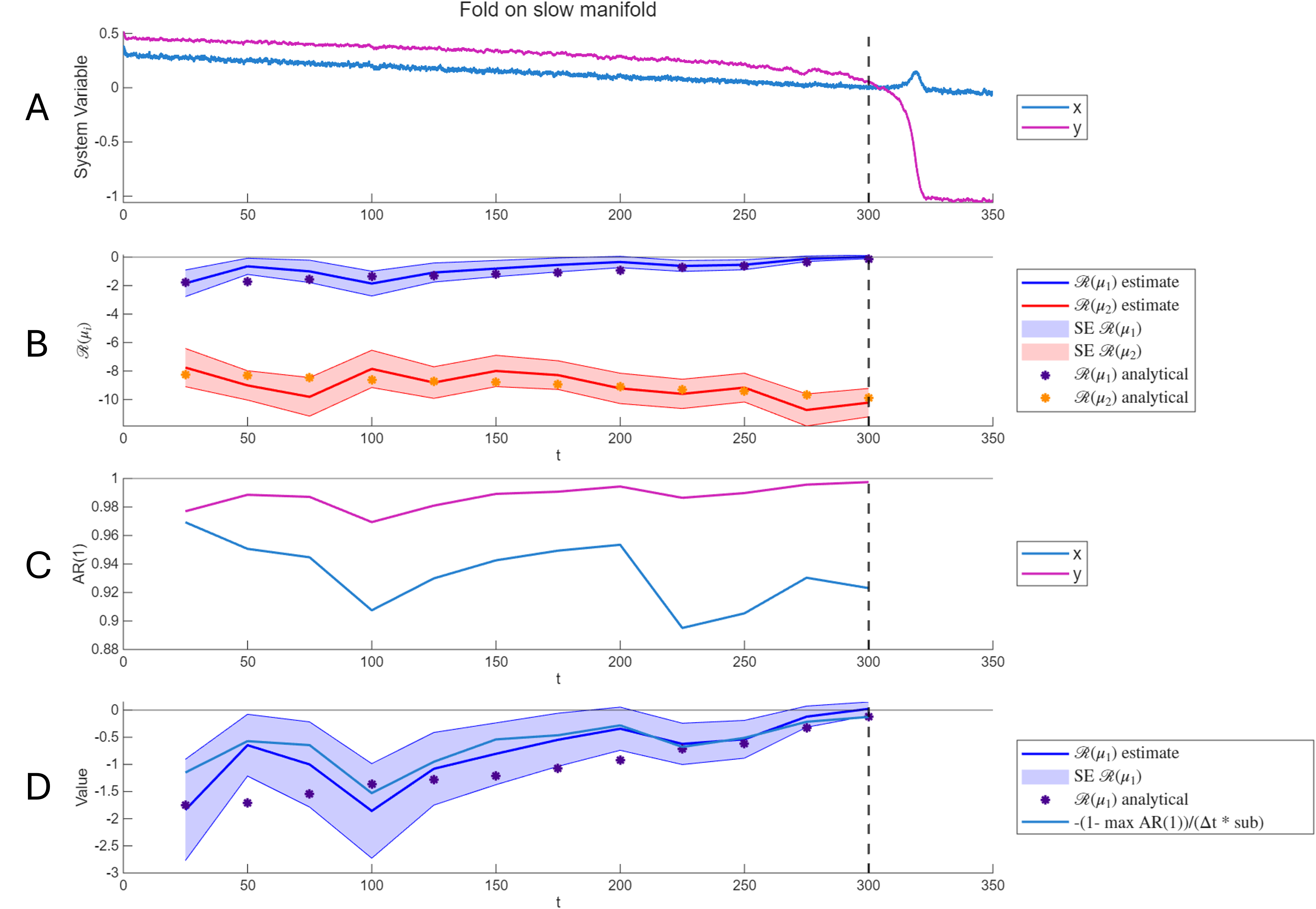}
\caption{Simulation for the fold system in (\ref{eq:fold_dxdy}) with parameters from Table \ref{tab:params}. In each plot, the black vertical line indicates the fold bifurcation point ($\lambda=0$). (A) Time series of $x$ (blue) and $y$ (magenta). (B) Real parts of VAR-estimated eigenvalues: thick blue and red lines show the leading and non-leading eigenvalues, with shaded regions indicating standard errors from Monte Carlo simulations. Stars denote analytical eigenvalues from the system Jacobian. (C) Lag-1 autocorrelation of $x$ and $y$. The same window length is used for both AR(1) and VAR calculations. (D) real part of the leading eigenvalue calculated using VAR, and AR(1) on $y$.}
\label{fig:foldevs}
\end{figure}
For the system (\ref{eq:fold_dxdy}), we simulate using the parameters in Table~\ref{tab:params}. This example features three timescales as $r\ll 1\ll \epsilon^{-1}$. Note that tipping occurs quickly in $y$ relative to the forcing, despite $y$ being considered the `slow variable'. As $t$ increases, $y$ decreases gradually before abruptly dropping when bifurcation-induced tipping occurs at $y=0$. Our study compares vector autoregression (VAR) with standard early warning signals, specifically lag-1 autocorrelation (AR(1)) for a univariate series. 

Figure \ref{fig:foldevs}(B) shows how the real parts of the eigenvalues ($\mathbb{R}(\mu_1,\mu_2)$) change as the system nears tipping. The Jacobian is calculated at each time $t$ at the end of the observed window. We determine the system's analytical eigenvalues and compare them with the VAR estimates. The analytical eigenvalues are plotted in purple and yellow for the leading and non-leading eigenvalues. This provides a clear visual comparison to assess whether the VAR model produces reasonable eigenvalues. 

The leading eigenvalue of the Jacobian $Df(x*)$ is the least stable, identified by having the smallest negative real part, denoted $-\alpha$ for $\alpha>0$. The parameter $\alpha$ measures the return rate: crossing zero signals a loss of stability. For the fold bifurcation, the VAR estimate for the leading eigenvalue is shown in blue. Over the series, the VAR estimate of the real part of the leading eigenvalue increases gradually as the system approaches the fold bifurcation, indicating a gradual loss of stability. The real part of the eigenvalue in red reflects the stable attraction to the slow manifold along the $x$-direction (Figure \ref{fig:foldevs}B). Both the VAR estimate and the analytical value for the non-leading eigenvalue decrease very gradually over the time series, suggesting a slight increase in stability as the system approaches the fold bifurcation, as verified by examining the eigenvalue structure, given by
\begin{equation}
    \mu_{1,2} = \frac{-1 \pm \sqrt{1-4\epsilon(2y+3y^2)}}{2\epsilon}.
\end{equation}
where, without loss of generality, we assume $\Re(\mu_1)\geq \Re(\mu_2)$.

Increasing the timescale separation (making $\epsilon$ smaller) between $x$ and $y$ would further decrease the real part of the non-leading eigenvalue, indicating greater attraction to the slow manifold in the fast direction. At the bifurcation point, the real part of the leading eigenvalue crosses zero. In each system we show, VAR estimation ends just before the tipping point, beyond this point the behaviour becomes nonlinear, which means the linear approximations are no longer informative. 

Figure \ref{fig:foldevs}(C) shows the lag-1 autocorrelation for the $x$ and $y$ time series. In AR(1) estimated from $y$ shows a non-monotonic rise appears, while the AR(1) estimated from $x$ does not show a trend indicating an approaching tipping point. Figure \ref{fig:foldevs}(D) directly compares the VAR estimate for the leading eigenvalue with a function of the lag-1 autocorrelation on $y$. In the presented two-dimensional fold bifurcation example with changes in the slow subsystem, the leading eigenvalue most accurately reflects future stability. Although the lag-1 autocorrelation on $y$ has a similar trend to the VAR's leading eigenvalue estimate, this similarity means one might argue that VAR provides no additional insight over AR(1) in this case. The comparison thus highlights that effective early warning with AR(1) depends on selecting the right observable; here, it is the $y$ variable.

\subsection{VAR for subcritical Hopf}

Analogously to the previous case, we consider the system in equation (\ref{eq:hopfz}) in Cartesian coordinate form as an SDE with additive noise:
\begin{equation}
	\begin{aligned}
		\label{eq:hopf_dxdy}
		dx &=\big[-\lambda x -\omega y +x^3 +xy^2 - x^5 -2x^3y^2 -xy^4\big]dt +\alpha_xdW_x\\
		dy&=\big[-\lambda y +\omega x +y^3 +x^2y -x^4y -y^5 -2x^2y^3\big]dt +\alpha_ydW_y\\
		d\lambda &=-rdt
	\end{aligned}
\end{equation}

\begin{figure}[ht]
    \centering
    \includegraphics[width=\textwidth]{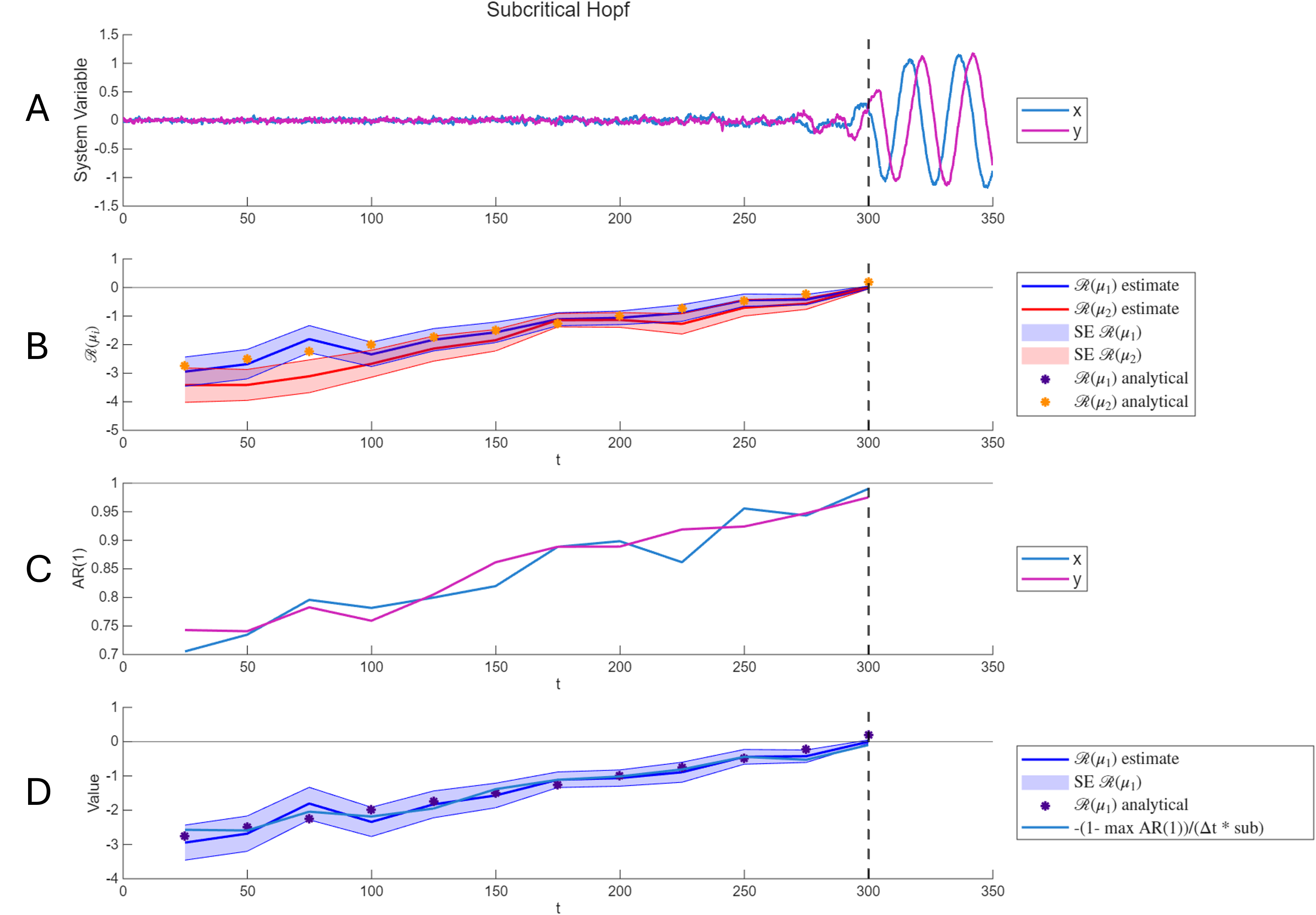}
    \caption{Simulation for the subcritical Hopf system in (\ref{eq:hopf_dxdy}) with parameters from Table \ref{tab:params}. In each plot, the black vertical line indicates the subcritical Hopf bifurcation point ($\lambda=0$). (A) Time series of $x$ (blue) and $y$ (magenta). (B) Real parts of VAR-estimated eigenvalues: thick blue and red lines show the leading and non-leading eigenvalues, with shaded regions indicating standard errors from Monte Carlo simulations. Stars denote analytical eigenvalues from the system Jacobian. (C) Lag-1 autocorrelation of $x$ and $y$. The same window length is used for both AR(1) and VAR calculations. (D) real part of the leading eigenvalue calculated using VAR, and AR(1) on $x$.}
    \label{fig:hopfevs}
\end{figure}

The numerical simulation of the subcritical Hopf bifurcation uses the parameters in Table~\ref{tab:params}. In this example, both $x$ and $y$ evolve on similar timescales. The time series in Figure \ref{fig:hopfevs}(A) shows both variables initially remain close to a stable equilibrium, with fluctuations around the origin increasing as the bifurcation is approached. Eventually, the system transitions to large-amplitude oscillations as the origin loses stability.

Figure \ref{fig:hopfevs}(B) presents the real parts of the VAR-estimated and analytical eigenvalues for the subcritical Hopf system. The VAR process typically estimates eigenvalues with equal real parts, reflecting a stable spiral equilibrium. At $t=75$, the VAR cannot resolve the imaginary component and instead produces two similar but distinct real parts. As time progresses and larger oscillations appear, the VAR model more reliably captures the imaginary component. When a complex conjugate pair is present, both eigenvalues are considered leading.

The ability to resolve the imaginary component improves when larger blocks and more frequent sampling are used, but at the expense of a longer time series. As $\lambda$ decreases, the system loses stability as the real part of the eigenvalue approaches zero when approaching the Hopf bifurcation at $t=300$. Once the system enters the regime of large-amplitude oscillations, the real part of the eigenvalue approaches zero, marking the transition to purely imaginary eigenvalues.

The AR(1) fit for both $x$ and $y$ in Figure \ref{fig:hopfevs}(C) also provides a clear indication of the loss of stability. This likely results from the block size used for AR(1) and VAR calculations being much larger than the oscillation period of the time series. For smaller block sizes, the imaginary component would be harder to resolve, and the increase in AR(1) would be less pronounced.

Figure \ref{fig:hopfevs}(D) compares the real part of the leading eigenvalue (either could be chosen in this case) with a function of the maximum AR(1). This panel makes clear that both methods consistently indicate a loss of stability as the subcritical Hopf bifurcation is approached.

\subsection{VAR for singular Hopf}

For the system (\ref{eq:singHopfxyslow}) with added noise, we consider the system of SDEs given by
\begin{equation}
	\begin{aligned}
		\label{eq:singhopf_dxdy}
		\epsilon dx &= \big[y-x^2(1+x)\big]dt + \alpha_xdW_x,\\
		dy &= \big[\lambda -x\big]dt +\alpha_ydW_y,\\
		d\lambda &=-rdt.
	\end{aligned}
\end{equation}

\begin{figure}[ht]
    \centering
    \includegraphics[width=\textwidth]{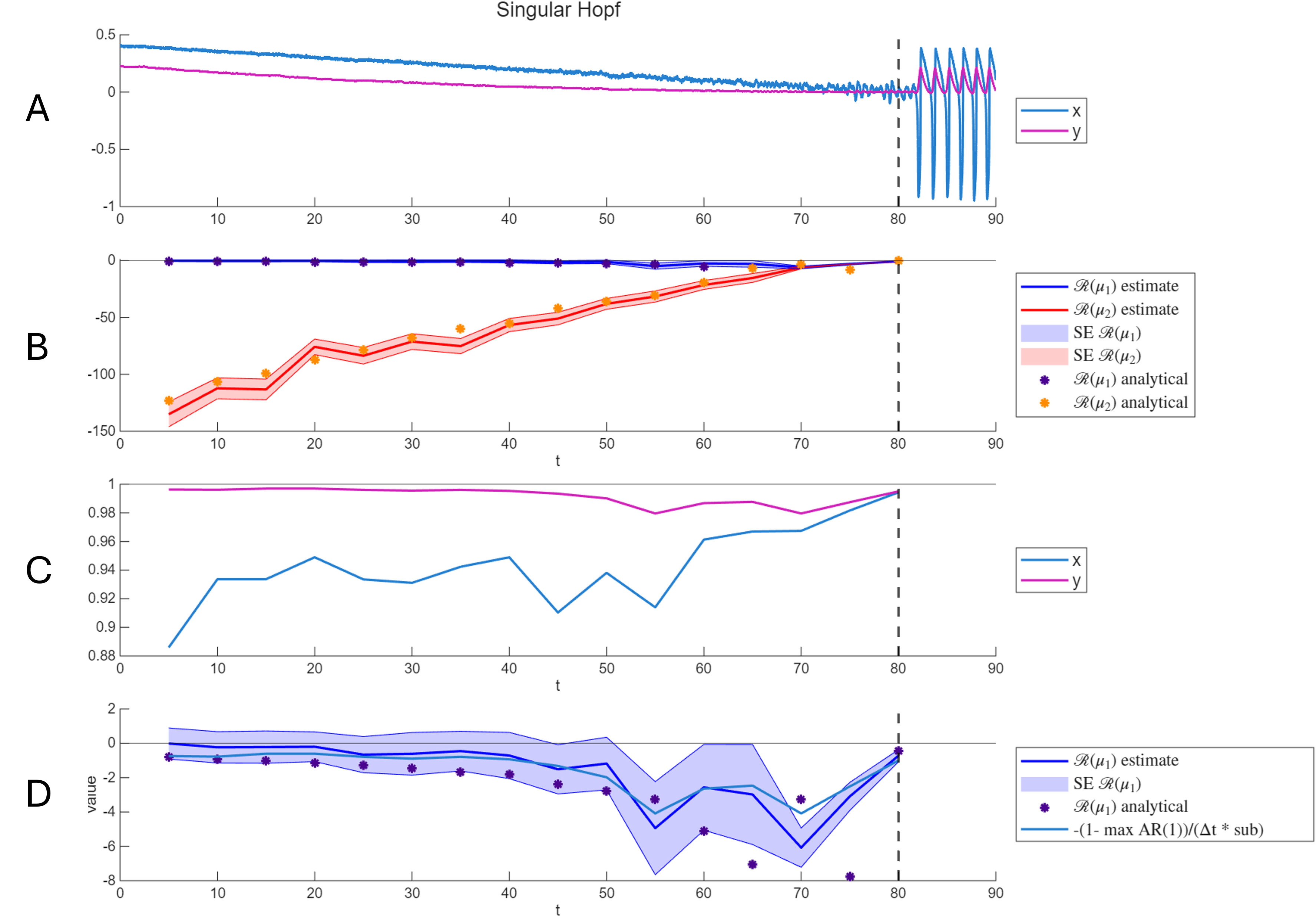}
    \caption{Simulation for the singular Hopf system in (\ref{eq:singhopf_dxdy}) with parameters from Table \ref{tab:params}. In each plot, the black vertical line indicates the singular Hopf bifurcation point ($\lambda=0$). (A) Time series of $x$ (blue) and $y$ (magenta). (B) Real parts of VAR-estimated eigenvalues: thick blue and red lines show the leading and non-leading eigenvalues, with shaded regions indicating standard errors from Monte Carlo simulations. Stars denote analytical eigenvalues from the system Jacobian. (C) Lag-1 autocorrelation of $x$ and $y$. The same window length is used for both AR(1) and VAR calculations. (D) Comparison of the real part of the leading eigenvalue calculated using VAR, and AR(1) on $x$.}
    \label{fig:singhopfevs}
\end{figure}

The implementation of the Singular Hopf bifurcation system uses parameter values given in Table \ref{tab:params} and we choose $\epsilon=0.01$ such that $r\ll 1\ll \epsilon^{-1}$ and so the fast dynamics and tipping are observed in the $x$ direction. Figure \ref{fig:singhopfevs}(A) illustrates a gradual decline in the $x$ value before a tipping point is reached, followed by a transition to an alternative state, specifically a relaxation oscillation. Similar behaviour is observed in the $y$ variable, but lesser in amplitude. For this specific set of parameters, the system visibly tips after the bifurcation point for $\lambda=0$ and $t=80$, indicated by the dashed black line. 

Figure \ref{fig:singhopfevs}(B) presents both the estimated and analytical eigenvalues for the singular Hopf system. The real part of the leading eigenvalue remains negative but steadily decreases in absolute value over the time series. A more apparent decrease is observed from approximately $t=40$ onward, which misleadingly suggests increasing stability, contrary to the actual decrease in stability as the singular Hopf bifurcation is approached. Due to the significant timescale separation, the non-leading eigenvalue initially has a substantially negative real part and steadily increases, eventually converging with the leading eigenvalue. The convergence of the real parts of the two eigenvalues marks the onset of oscillatory behaviour, resulting in a complex conjugate eigenvalue pair. This real part increases slowly until it crosses zero at the tipping point, which occurs after the singular Hopf bifurcation is reached.

Figure \ref{fig:singhopfevs}(C) shows AR(1) for $x$ and $y$ variables. AR(1) for $y$ decreases over the time series, while AR(1) for $x$ increases non-monotonically. Figure \ref{fig:singhopfevs}(D) compares the leading eigenvalue with a function of the maximum AR(1). On this scale, both decreasing AR(1) and the VAR estimate of the real part of the leading eigenvalue suggest increasing stability, failing to indicate the approaching bifurcation.

For the singular Hopf bifurcation, it is the non-leading eigenvalue (red in Figure \ref{fig:singhopfevs}(B)) that provides insight into the system’s future stability. Examining only the AR(1) of the univariate $x$ time series—despite $x$ being the 'tipping' variable—does not reveal the clearest trend. This example demonstrates the value of looking beyond AR(1) and univariate time series. Although the non-leading eigenvalue does not evolve linearly, one can apply linear regression to its real part to estimate the tipping time. Once the real parts of the eigenvalues have converged, the rate at which their real parts increase slows, so applying linear regression before this point will likely predict a tipping time that precedes the actual event. Nevertheless, this approach is still valuable as a starting point for prediction in a system with a large timescale separation.

\section{Discussion}
\label{sec:discuss}

For a slowly changing system with low-amplitude stochastic forcing, we can aim to extract indicators of stability by fitting to various autoregressive models such as vector autoregressive $\VAR(k,1)$ models considered here or, for example, autoregressive moving average (ARIMA) models \cite{rodal2022dynamical}. The simplest option AR(1) (corresponding to $\VAR(1,1)$) assumes one-dimensional dynamics and so, while useful for fold bifurcations, is not able to distinguish between the different types of transition we discuss here. We suggest that examining $\VAR(k,1)$ will be informative for this, especially if we have access to $k$ independent time series from the system we wish to understand. Even if this is not the case and we have only one time series, it should be possible to extract additional eigenvalues by fitting to $\VAR(1,k)$ for $k>1$, which is presumably equivalent to performing a Takens embedding and the $\VAR(k,1)$ on this.

In this study, we have highlighted three distinct types of generic codimension-1 instabilities with multiple timescales: a fold bifurcation on the slow manifold, a subcritical Hopf bifurcation, and a singular Hopf bifurcation. Many studies of early warnings of bifurcation-induced tipping points limit themselves to fold bifurcations. This means they miss the possibility that the bifurcation associated with the tipping point may be of Hopf type. Moreover, in a multiscale system, a singular Hopf bifurcation may occur generically.  We have chosen these in particular because they are common mechanisms for critical transitions and are therefore important in tipping point theory. We have shown that, for slowly ramped systems with a small-amplitude additive noise, it is possible to identify the leading and nonleading eigenvalues of the underlying system using a VAR model.

For the fold bifurcation on the slow manifold, the leading eigenvalue is the best indicator of when the system will cross the instability threshold. In this case, the VAR estimate of the leading eigenvalue exhibits a very similar trend to the univariate AR(1) model for $y$. For the subcritical Hopf bifurcation, we demonstrated the behaviour of the eigenvalues that would appear in a system with a timescale separation, the only difference being that there would be an additional eigenvalue on the order of $O(\epsilon^{-1})$, that would be consistently negative, similar to the fold case. The VAR estimation identifies two eigenvalues with the same real component and therefore resolves the imaginary component associated with oscillatory behaviour near the equilibrium preceding the subcritical Hopf bifurcation. Because the block size is much larger than the period of oscillation, the AR(1) strongly indicates a loss of stability as the system approaches the bifurcation. The usefulness of VAR is most clearly demonstrated in the singular Hopf bifurcation example. The leading eigenvalue is only useful as an indicator of tipping very close to the bifurcation. From far away, the nonleading eigenvalue indicates decreasing stability. The AR(1) applied to the $x$ and $y$ time series in this case would not provide a strong indication of tipping if this were the only method we used.

Since a VAR model generalises an AR(1) model to multiple variables, it relies on analogous assumptions, extended to the multivariate setting, in order to detect early warning signals. These include and are not limited to: requirement for the system to be in steady state, changes in the system are due to a slowly varying parameter and that the parameter can be considered stationary with respect to the other timescales of the system dynamics, noise amplitide is small, that the system can be modelled by an ODE for the deterministic part and a Gaussian noise process for the remainder, and that chnages in the system are attributable to a bifurcation rather than noise induced tipping for example \cite{Williamson_2015}. In addition, estimating the linearization using VAR should help avoid problems such as those outlined in \cite{Morr2024}, where ``internal noise interference'' can mask AR(1) trends for some observables.

The $\VAR(k,1)$ method faces many of the same challenges as AR(1), in particular, choice of detrending method, sampling times and length of the sliding window. Although sensitivity analyses can help assess the impact of these factors on observed early warning signals, issues such as subsampling, sparse data sets, and insufficiently long time series remain significant hurdles. The method requires a long time series, so it may be ineffective for tipping in systems with short time horizons. Moreover, as VAR estimates more parameters from the same dataset, a greater amount of data is needed to achieve comparable levels of uncertainty for each estimate, further increasing the demand for data.

One challenge is that $\VAR(k,1)$ and AR(1) may interpret slow changes in the system as additional stable eigenvalues. Hence, it is important to provide an estimate of the uncertainty in the eigenvalues. Other approaches to estimate eigenvalues include, for example, Koopman/transfer operator/response theory methods (eg \cite{lucarini2026general,zagli2026bridging}) though in many cases it is not easy to estimate the uncertainty in the spectrum. We suggest that the estimated eigenvectors will give further valuable information near instability, though we do not explore this further here.

Transforming the estimated eigenvalues into an early warning signal requires further steps, firstly to estimate a trend (with uncertainty) and then to extrapolate this to a crossing of the instability threshold. If we know the system is in a regime where a bifurcation normal form is valid and the forcing is approximately linear, the extrapolated trend can yield a skilful classifier for tipping within some time horizon. This is considered in \cite{Ashwin2025} for the fold bifurcation and should be possible for subcritical Hopf, and may be possible for the singular Hopf: as noted in \cite{hobden2025regularization}, the leading eigenvalue is not informative, but we surmise that fitting to a slow-fast normal form may be informative.

There is one additional generic case in multiscale systems that we do not consider here. This is a Fast Hopf bifurcation with an additional dimension with stable slow dynamics. This occurs in at least three dimensions when a complex pair of fast eigenvalues crosses the imaginary axis at a subcritical Hopf bifurcation. This case is interesting because the complex pair of eigenvalues only becomes leading very close to the bifurcation. Put another way, similar to the singular Hopf, the leading eigenvalue is not informative of the trend except very near the transition. We do not consider bifurcations of periodic or chaotic attractors, but note that in these cases the instantaneous Jacobian is not necessarily informative for the stability - longer term averages will need to be integrated over trajectories for this. For more complex attractors, other transitions can occur but we do not consider these here \cite{Alkhayoun2019,Alkhayuon2021}. 

In summary, any application of this method to provide skilful early warnings of tipping points for time series ``in the wild'' will face many challenges. These include choosing high signal-to-noise observations, estimating eigenvalues tightly enough to measure trends with good reliability, and understanding when and if these trends can be extrapolated. However, most of these challenges are already faced by AR(1) approaches, while $\VAR(k,1)$ methods can potentially overcome some of the structural limitations of AR(1).

\subsection*{Data availability statement}

The code required to create the figures in this paper is downloadable from \url{https://github.com/Bryony-Hobden/VAR_EWS}.

\subsection*{Acknowledgements}

We thank Mark Williamson and George Datseris for interesting conversations and useful feedback. BH thanks the Faculty of Environment, Science and Economy  of the University of Exeter for support via a PhD studentship. PA and PR thank ClimTip and AdvanTip for funding support/. The ClimTip project has received funding from the European Union's Horizon Europe research and innovation programme under grant agreement No. 101137601: Funded by the European Union. Views and opinions expressed are however those of the author(s) only and do not necessarily reflect those of the European Union or the European Climate, Infrastructure and Environment Executive Agency (CINEA). Neither the European Union nor the granting authority can be held responsible for them. The Advanced Research and Invention Agency (ARIA) fund the AdvanTip project (grant no. SCOP-PR01-P003).

For the purpose of open access, we have applied a Creative Commons Attribution (CC BY) license to any Author Accepted Manuscript version arising from this submission.

\bibliographystyle{plain}
\bibliography{EWSrefs}

\newpage
\appendix

\section{Delta method for propagating standard errors}
\label{app:delta}

We briefly outline the ``delta'' or first-order error approximation from \cite[Section 5.5.4]{Casella_Berger_2024}.
Let $\mathbf{X}$ be a vector of random variables $\mathbf{X} =(X_1,\dots X_n)$ with mean $\mathbf{X}^*=(X_1^*,\dots,X_n^*)$. Suppose there is a differentiable function $F(\mathbf{X})$ for which we want to estimate mean and variance. Let
\begin{equation}
F'_i(\mathbf{X}^*) = \frac{\partial}{\partial X_i} F|_{\mathbf{X}=\mathbf{X}^*}.
\end{equation}
The first-order Taylor expansion of $F$ around $\mathbf{X}^*$ is
\begin{equation}
    F(\mathbf{X}) = F(\mathbf{X}^*) + \sum_{i=1}^n F'_i(\mathbf{X}^*)(X_i-X^*_i) + O(|\mathbf{X}-\mathbf{X}^*|^2)
\end{equation}
If we assume higher-order terms to be less significant, we have the approximation
\begin{equation}
    F(\mathbf{X}) \approx F(\mathbf{X}^*) + \sum_{i=1}^n F'_i(\mathbf{X}^*)(X_i-X^*_i).
    \label{eqn:DeltaTaylor}
\end{equation}
This means that the mean of $F(\mathbf{X})$ is simply $F(\mathbf{X}^*)$. We can estimate the variance of $F(\mathbf{X})$ using equation \ref{eqn:DeltaTaylor} by
\begin{align}
\label{eqn:Varf}
     \var(F(\mathbf{X})) &\approx \mathbb{E}\big([F(\mathbf{X})-F(\mathbf{X}^*)]^2\big)\\
     &\approx \mathbb{E}\Big(\sum_{i=1}^n F'_i(\mathbf{X}^*)(X_i-X^*_i)\Big)\\
     &= \sum_{i=1}^n (F'_i(\mathbf{X}^*))^2 \var(\mathbf{X}) + 2\sum_{i>j}F'_i(\mathbf{X}^*)F'_j(\mathbf{X}^*)\cov(X_i,X_j).
\end{align}

In order for the delta method to give a good estimate, the function $F$ must be smooth and differentiable near the mean. Additionally, the variables should be independent (no covariance) or the covariance structure must be known. In our case, we assume there is no covariance and so equation \ref{eqn:Varf} simplifies to
\begin{align}
    \var(F(\mathbf{X})) &\approx \sum_{i=1}^n (F'_i(\mathbf{X}^*))^2 \var(\mathbf{X})\\
    &\approx \sum_{i=1}^n (F'_i(\mathbf{X}^*))^2 SE(\mathbf{X}^2).
\end{align}

\end{document}